%% file: main.tex
\documentclass[10pt, reqno]{amsart}

\usepackage{geometry}
\usepackage{url}
\usepackage{amssymb}
\usepackage[english]{babel}
\usepackage{colortbl}
\usepackage{enumerate}
\usepackage{graphicx} 
\usepackage{float}
\usepackage{mathrsfs}
\usepackage{cite}
\usepackage{kpfonts}
\usepackage{tcolorbox} 
\usepackage{enumitem}
\usepackage{bbm}
\usepackage{amsmath}
\usepackage{amsfonts}
\usepackage{latexsym}
\usepackage{mathtools}
\usepackage{rotating}     
\usepackage{pdflscape}
\usepackage{siunitx}
\usepackage{diagbox}

\allowdisplaybreaks

\usepackage[colorlinks=true,urlcolor=green!50!black, citecolor=red!60!black,linkcolor=blue!75!black,linktocpage,pdfpagelabels,bookmarksnumbered,bookmarksopen]{hyperref}
\usepackage[capitalize, noabbrev]{cleveref}

\graphicspath{{Figures}}
\usepackage{float}

\usepackage{tikz}
\usetikzlibrary{arrows, automata, shapes, calc, patterns}

\usepackage{pgf}
\usepackage{pgfplots}
\pgfplotsset{compat=1.18}

\usepgfplotslibrary{external}
\usepackage{adjustbox}

\usepackage{etoolbox}

\usepackage{mathtools}
\DeclarePairedDelimiter\abs{\lvert}{\rvert}
\DeclarePairedDelimiter\norm{\lVert}{\rVert}
\DeclarePairedDelimiterX{\inner}[2]{\langle}{\rangle}{#1, #2}
\DeclarePairedDelimiter\floor{\lfloor}{\rfloor}

\newcommand{\R}{\mathbb{R}}

\newcommand{\Lp}[1]{\mathit{L}^{#1}}

\newcommand{\normLp}[2]{\norm{#1}_{\Lp{#2}}}

\renewcommand{\uplus}{u^+}
\newcommand{\uminus}{u^-}
\newcommand{\upm}{u^{\pm}}
\newcommand{\upxt}{u^+(x,t)}
\newcommand{\umxt}{u^-(x,t)}
\newcommand{\uxt}{u(x,t)}
\newcommand{\lambdaplus}{\lambda^+}
\newcommand{\lambdaminus}{\lambda^-}
\newcommand{\lambdapm}{\lambda^{\pm}}
\newcommand{\ypm}{y^{\pm}}

\newtheorem{theorem}[]{Theorem}[section]

\numberwithin{equation}{section}

\newtheorem{definition}[theorem]{Definition}

\newtheorem{remark}[theorem]{Remark}

\begin{document}

\input{Title_Abstract}

\input{Introduction}

\input{Hyperbolic_model}

\input{Numerical_scheme}

\input{Choice_numerical_scheme}

\input{Conclusion}

\input{Appendix}

\section*{Acknowledgments} The authors (TTL and RE) acknowledge funding from a Region Bourgogne Franche-Comt\'{e} ``Accueil de Nouvelle \'{E}quipe de Recherche (ANER) 2022" grant number   \verb|FC22070.LMB.CL.|, and funding through a University of Franche-Comt\'{e} "Chrysalide: soutien aux nouveaux arrivants" (for RE).
Computations have been performed on the supercomputer facilities of the Mésocentre de calcul de Franche-Comté and on the computation server of Laboratoire de Math\'{e}matiques de Besan\c{c}on. The authors would like to thank Julien Yves ROLLAND for his support on the computation servers.

\renewcommand{\bibname}{References}
\bibliography{References_Postdoc}
\bibliographystyle{plain}

\end{document}

%% file: Title_Abstract.tex
\title[]{Numerical challenges for the understanding of localised solutions with different symmetries in non-local hyperbolic systems}

\author[T.T. Le and R. Eftimie]{Thanh Trung Le and Raluca Eftimie}

\address{
	\vspace{-0.25cm}
	\newline
	\textbf{{\small Thanh Trung Le}} 
        \newline \indent Universit\'{e} de Franche-Comt\'{e}, CNRS, LmB, F-25000 Besan\c{c}on, France}
\email{thanh\_trung.le@univ-fcomte.fr} 

\address{
	\vspace{-0.25cm}
	\newline
	\textbf{{\small Raluca Eftimie}} 
        \newline \indent Universit\'{e} de Franche-Comt\'{e}, CNRS, LmB, F-25000 Besan\c{c}on, France}
\email{raluca.eftimie@univ-fcomte.fr}

\subjclass[2020]{35, 37, 45, 65, 68, 92}
\keywords{Snake-and-ladder bifurcation; Nonlocal hyperbolic systems; Ecological aggregations; Finite volumes schemes; Convergence and non-convergence of numerical schemes; Numerical transient solution; Numerical steady-state solution.}
\maketitle

\begin{abstract} 
    We consider a one-dimensional nonlocal hyperbolic model introduced to describe the formation and movement of self-organizing collectives of animals in homogeneous 1D environments. Previous research has shown that this model exhibits a large number of complex spatial and spatiotemporal aggregation patterns, as evidenced by numerical simulations and weakly nonlinear analysis. In this study, we focus on a particular type of localised patterns with odd/even/no symmetries (which are usually part of snaking solution branches with different symmetries that form complex bifurcation structures called snake-and-ladder bifurcations).
    To numerically investigate the bifurcating solution branches (to eventually construct the full bifurcating structures), we first need to understand the numerical issues that could appear when using different numerical schemes. To this end, in this study, we consider ten different numerical schemes (the upwind scheme, the MacCormack scheme, the Fractional-Step method, and the Quasi-Steady Wave-Propagation algorithm, combining them with high-resolution methods), while paying attention to the preservation of the solution symmetries with all these schemes. We show several numerical issues: first, we observe the presence of two distinct types of numerical solutions (with different symmetries) that exhibit very small errors, which might initially suggest that we have reached a steady-state solution, but this is not the case (this also implies an extremely slow convergence);  second, in some cases, none of the investigated numerical schemes converge, posing a challenge for the development of numerical continuation algorithms for nonlocal hyperbolic systems; lastly, the choice of the numerical schemes, as well as their corresponding parameters such as time-space steps, exert a significant influence on the type and symmetry of bifurcating solutions. To conclude we emphasize that if we want to construct numerically bifurcation diagrams for these localised solutions with different symmetries, the resulting bifurcations may vary when different numerical schemes and/or corresponding parameters are employed.
\end{abstract}

%% file: Introduction.tex
\section{Introduction}

The study of animal aggregations has been investigated intensively over the past fifty years \cite{Okubo-Grunbaum-Keshet-2001, Flierl-Grunbaum-Levin-Olson-1999, Gueron-Levin-Rubenstein-1996, Mogilner-Keshet-1999, Reynolds-1987}. One of the most studied aspects of these aggregations is the spatial and spatiotemporal patterns exhibited by them:  zigzagging flocks of birds \cite{Feder-2007}, or milling schools of fish, for instance, \cite{Parrish-Keshet-1999} and \cite{Lukeman-Li-Keshet-2009}. To investigate the biological mechanisms necessary for the formation and persistence of these patterns, a wide range of mathematical models have been proposed. There are two main classes of mathematical models used for animal aggregations: 1) individual-based models (Lagrangian approach, microscopic models), which track the movements of all individuals in the group, and 2) partial differential equations (PDE) models, formulated as evolution equations for the population density ﬁeld. Despite the complex group patterns displayed by the individual-based models (e.g., swarms, tori, polarized groups, see \cite{Couzin-Krause-James-Ruxtion-2002} and the reference therein), the lack of techniques to investigate them causes difﬁculties in understanding some of these patterns, see \cite{Eftimie-Vries-Lewis-2009}. Hence, the main method for this approach focuses on numerical simulations with the purpose of comparing the simulated aggregation patterns with the available data, and thus ascertaining the correctness of the micro-scale level assumptions incorporated into the models, for example, see \cite{Aldana-Dossetti-Huepe-Kenkre-Larralde-2007, Czirok-Barabasi-Vicsek-1999, Czirok-Stanley-Vicsek-1997, Vicsek-Czirok-Jacob-Cohen-Shochet-1995}. The second approach can be classified into two categories: kinetic models (mesoscopic models; see \cite{Othmer-Dunbar-Alt-1988, Fetecau-2011}, also see \cite{Eftimie-2012}) and continuum models (Eulerian approach, macroscopic models) which is better represented in the literature on animal aggregations, with a diverse range of parabolic models \cite{Mogilner-Keshet-1999, Topaz-Bertozzi-Lewis-2006} and hyperbolic models \cite{Eftimie-Vries-Lewis-Lutscher-2007, Eftimie-Vries-Lewis-2007, Lutscher-2002, Pfistner-1990}. For continuum models, the main method focuses on the analytical results, e.g., showing the existence of particular solutions exhibited by these models, showing the existence of different types of bifurcations that give rise to different solution types, or trying to connect biological interactions at the micro-scale and macro-scale levels (through various homogenisation approaches). 
There are, however, models that combine analytical results with numerical simulations, for example, see \cite{Chuang-Orsogna-Marthaler-Bertozzi-Chayes-2007, Fetecau-Eftimie-2010, Buono-Eftimie-2014-SIAM, Eftimie-Vries-Lewis-2009}. \\ 

In \cite{Eftimie-Vries-Lewis-2007} and \cite{Eftimie-Vries-Lewis-Lutscher-2007}, the authors have proposed a nonlocal hyperbolic model that introduces a general framework to incorporate different communication mechanisms to study the formation of animal groups. In particular, these communication mechanisms inﬂuence the social interactions between individuals, namely attraction towards other members of the group that are far away, repulsion from those that are nearby, and a tendency to align with those neighbors that are at intermediate distances. The resulting model, which actually comprises many submodels, is very rich in spatial and spatiotemporal patterns. Numerical simulations have shown at least 10 different patterns, including stationary and traveling pulses, traveling trains, zigzag pulses, breathers, ripples, and feathers. Two particular patterns, traveling trains, and stationary pulses, were shown to arise through subcritical Hopf and Steady-state (codimension-one) bifurcations  \cite{Eftimie-Vries-Lewis-2009}. Other patterns, such as ripples, were shown to arise through codimension-two bifurcations \cite{Buono-Eftimie-2014-MMMA, Buono-Eftimie-2014-SIAM}. In this study, we focus on a particular type of localised patterns with odd/even/no symmetries (which are usually part of snaking solution branches with different symmetries that form complex bifurcation structures called snake-and-ladder bifurcations). The snake-and-ladder bifurcations have been observed in various simple fluid dynamics models such as Swift-Hohenberg equation (see \cite{Avitabile-Lloyd-Burke-Knobloch-Sandstede-2010, Beck-Knobloch-Lloyd-Sandstede-Wagenknecht-2009, Burke-Knobloch-2006, Burke-Knobloch-2007-Chaos, Burke-Knobloch-2007-PLA, Liu-Xu-2017} and the reference therein), discrete bistable Allen-Cahn equation \cite{Taylor-Dawes-2010} and more recently in simple integro-differential equations \cite{Schmidt-Avitabile-2020}. Normally, the snake-and-ladder bifurcation is structured by two intertwined snaking branches corresponding to even and odd localised solutions. These branches include both stable and unstable branches. Near each fold, the localised patterns undergo a pitchfork bifurcation at which a branch of asymmetric solutions emerges that connects the two snaking branches. The ladder branches may be unstable or stable. For example, we refer to \cite[Figure 1]{Avitabile-Lloyd-Burke-Knobloch-Sandstede-2010}.\\

Understanding the dynamics of a system describing a physical, biological, or engineering problem, where changes in parameter values can lead to changes in system dynamics, is an important fundamental aspect of applied mathematics. In addition to using analytical approaches to investigate the impact of parameter changes on the number and stability of different model states, numerical simulations are being used to visualize the dynamics of the model as one or multiple parameters are varied. However, such numerical simulations could be very time-consuming due to the manual inspection of solution trajectories to detect transient and asymptotic dynamics, multiple stable steady states, bifurcation points, stable and unstable manifolds, etc. 
As an alternative, numerical continuation algorithms are being used to automatically detect invariant sets, bifurcations, metastable states, etc. Numerical continuation algorithms for ordinary differential equations (ODEs) are now an established tool for bifurcation analysis in dynamical systems \cite{Kuehn-2015}. For ODEs, there are several standard good software packages (e.g., XPP \cite{XXP}, AUTO-07p \cite{AUTO-07p}, MatCont \cite{MATCONT}, PyDSTool \cite{PyDSTool}, etc). However, for partial differential equations (PDEs) there are no standard numerical continuation toolboxes that would cover a broad range of different classes of PDEs automatically \cite{Kuehn-2015}. The toolboxes that exist (e.g., \textit{pde2path} package in \textsc{MATLAB}) are usually developed for specific classes of partial differential equations (mainly of elliptic or parabolic types), using numerical methods most suitable for those classes. The simplest approach, but also very expensive computationally, is to implement a discretization scheme most suitable for a class of equations, and then add a predictor-corrector continuation algorithm on top of the discretization scheme. However, changing parameters could lead to changing the type and characteristics of the PDEs, which might require also changes in the numerical method used to track the solution branches. One way of solving this particular problem is to glue together different numerical methods developed for different classes of equations \cite{Kuehn-2015}.

Most of these numerical bifurcation studies have been applied to local PDEs of parabolic and elliptic types (mainly arising in fluid dynamics), in which the bifurcations diagrams are obtained with the help of standard continuation software (e.g., XPP, AUTO-07p, MatCont, PyDSTool, etc). More recently, they have started to be applied also to new mathematical models describing various biological and medical problems. Some of these models are also of non-local type, describing long-distance interactions between the various components of the system  \cite{Schmidt-Avitabile-2020}. However, the field of mathematical approaches to biology/medicine is very vast, with new complex nonlinear and non-local mathematical models (sometimes described by PDEs of hyperbolic type) that exhibit interesting patterns and bifurcations being developed continuously. Therefore, new numerical approaches need to be developed to investigate the bifurcations of these patterns in these new PDEs (see for example, the Galerkin schemes and Fourier collocation scheme introduced in \cite{Schmidt-Avitabile-2020}). This is particularly important for models of hyperbolic type, where the center manifold theorem (required for bifurcation theory) does not always hold true. \\

\begin{figure}[!ht]
    \centering
    \resizebox{0.9 \textwidth}{!}{\input{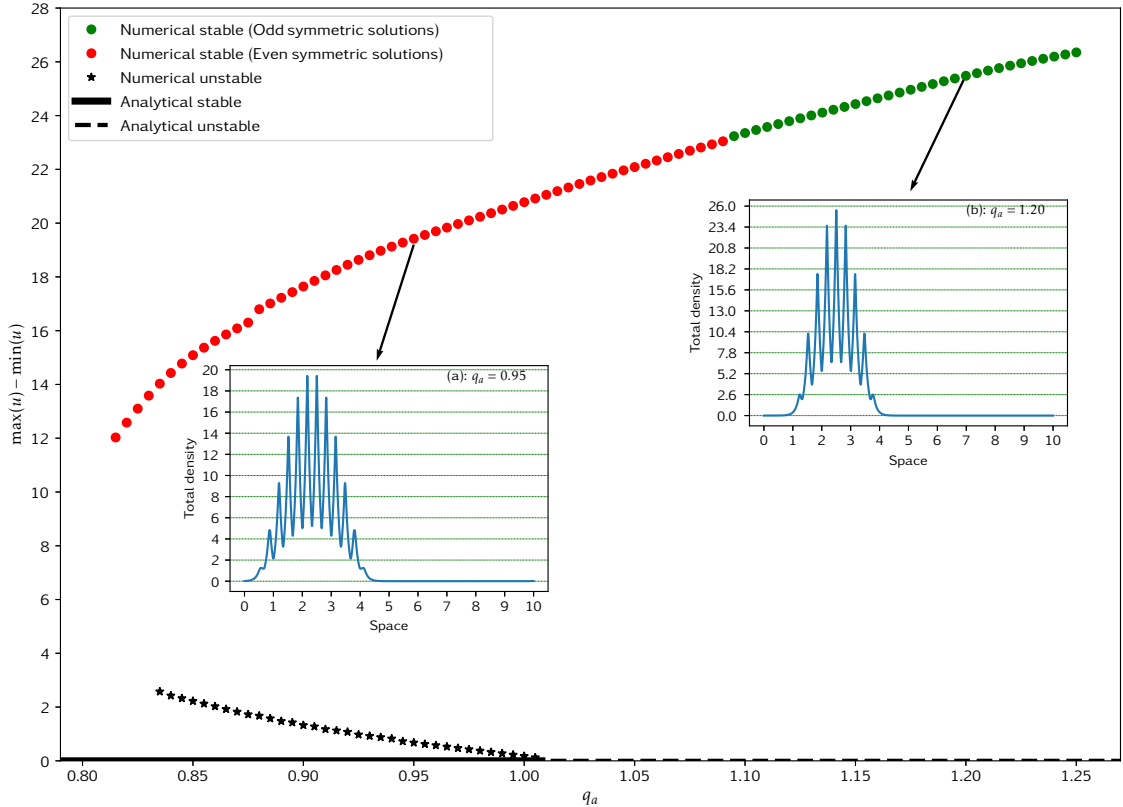}}
    \caption{Numerical bifurcation diagram obtained as we vary the magnitude of attraction $q_{a}\in [0.83,1.25]$. For $q_{a}<1.008$ the zero-amplitude spatially-homogeneous solution is stable, while for $q_{a}>1.008$ this zero-amplitude solution is unstable. At $q_{a}=1.01$ an unstable solution bifurcates sub-critically towards lower $q_{a}$ values. The starred curve bifurcating to the left is a separatrix curve acting as a boundary between the lower region where perturbations of spatially-homogeneous steady state with amplitudes below this curve are decaying towards zero, and perturbations with amplitudes above this value grow into the high-amplitude heterogeneous patterns shown on the upper branch (the filled circles). The upper branch shows the stable high-amplitude even symmetric solutions (for $q_{a}<1.095$) and the stable high-amplitude odd symmetric solutions (for $q_{a} \geq 1.095$).}
    \label{fig:Fi3_WeaklyNonlinear}
\end{figure}
In this paper, we focus on numerical approaches to simulate the localised patterns with odd/even/no symmetry (that likely form the branches of snake-and-ladder bifurcations) in the nonlocal hyperbolic systems for ecological aggregations introduced in \cite{Eftimie-Vries-Lewis-2007} and  \cite{Eftimie-Vries-Lewis-Lutscher-2007}. Specifically, we will consider here only one of the submodels introduced in \cite{Eftimie-Vries-Lewis-2007}, where information received from all neighbors are taken into account for attractive and repulsive interactions, while for alignment, only information from neighbors moving toward an individual is considered. But for simplicity, we assume that attraction and repulsion are the only possible social interactions (i.e., the magnitude of alignment is zero). 
Previous research \cite{Eftimie-Vries-Lewis-2007} has shown that this particular model exhibits stationary localised patterns, characterised by different symmetries. Using a weakly nonlinear analysis combined with numerical simulations, in \cite{Eftimie-Vries-Lewis-2009} the authors showed that the stationary pulses arise through a real subcritical (i.e., unstable) bifurcation from the spatially homogeneous steady state. In \cref{fig:Fi3_WeaklyNonlinear} we show the changes in these steady-state homogeneous solutions (with zero amplitudes) and the steady-state heterogeneous solutions (with non-zero amplitudes) as we vary the magnitude of attraction $q_{a}$. We observe here two different types of stable localised solutions on the high-amplitude branch: a stable solution with even symmetry (red circles, for $q_{a}<1.095$) and a stable solution with odd symmetry (green circles, for $q_{a}>1.0.95$). We also note two different types of unstable localised solutions (with odd and even symmetries) on the low-amplitude branch that bifurcates sub-critically at $q_{a}=1.01$.

To numerically understand these localised patterns with odd/even/no symmetries (which are probably part of a snake-and-ladder bifurcation structure, formed of two ``snaking" solution branches with odd and even symmetries connected by ``ladder" branches with no symmetry), we employ different numerical schemes. Here we focus on the first-order upwind scheme (see, e.g., \cite[Chapter 4]{LeVeque-2002-Finite-volume}), the second-order MacCormack scheme (see \cite{MacCormack-2003}, see also \cite{Helbing-Treiber-1999}), or the Fractional-Step method (see, e.g., \cite[Chapter 17]{LeVeque-2002-Finite-volume}). Our goal is to find a spatially-heterogeneous steady-state solution. However, a numerical challenge arises when approaching the steady state, as the flux gradient becomes almost balanced by the source term. As a result, the magnitude of the time derivative of the solution becomes much smaller than that of the spatial derivative of the flux function which is comparable to the magnitude of the source term. This implies that many numerical methods, such as the Fractional-Step method, have difﬁculty preserving such steady states and cannot accurately capture small perturbations around them. To overcome that, we consider the Quasi-Steady Wave-Propagation Algorithm (QSA) introduced in \cite{LeVeque-1998} and combine the algorithm with high-resolution methods. One of the new methods for balance laws (called the f-wave method) introduced by \cite{Bale-Leveque-Mitran-Rossmanith-2002}, is efficient with spatially varying flux functions and/or spatially varying source terms. For our problem, it is equivalent to QSA.

In this study, we investigate the impact of different numerical schemes on the type of solutions (even symmetric, odd symmetric, or non-symmetric), which could lead to changes in snake-and-ladder bifurcation diagrams. (Note that here we focus only on the choice of the numerical schemes; the construction of the bifurcation diagrams will be presented in a future paper). We will show here several numerical issues in our study. First, we observe the presence of two distinct types of numerical solutions that exhibit very small errors between subsequent time steps, which might initially suggest that we have reached a steady-state solution, but this is not the case. This also implies an extremely slow convergence.  Second, in some cases, none of the investigated numerical schemes converge, posing a challenge to the numerical analysis. Lastly, we have discovered that the choice of numerical schemes, as well as their corresponding parameters such as time-space step and initial conditions, exert a significant influence on the type and symmetry of bifurcating solutions. As a result, the resulting bifurcation diagrams may vary when different numerical schemes and/or corresponding parameters are employed.

The paper is structured as follows. In \cref{sec:model}, we provide a brief overview of the nonlocal hyperbolic model under analysis,  the spatially homogeneous steady states, and the localised heterogeneous solution. \cref{sec:numerical_scheme} introduces the numerical schemes that have been investigated in this study. \cref{sec:choice_scheme} represents the main focus of the paper, where we discuss the numerical challenges encountered. In \cref{sec:conclusion}, we summarize our findings and draw conclusions, along with highlighting some open problems for further research. The Appendix includes figures that enhance the clarity of the issues discussed in \cref{sec:choice_scheme}.

%% file: Hyperbolic_model.tex
\section{Nonlocal hyperbolic model and Localized spatially-heterogeneous solutions with different symmetries} \label{sec:model}

\textbf{Nonlocal hyperbolic model.} Following \cite{Eftimie-Vries-Lewis-2007} and \cite{Eftimie-Vries-Lewis-Lutscher-2007}, we consider the following one-dimensional hyperbolic model to describe the evolution of densities of right-moving ($\uplus$) and left-moving ($\uminus$) individuals in a generic ecological population on a 1D spatial domain (i.e., a domain much longer than wide):
\begin{align}
 \begin{cases}\label{eqn:main}
		\partial_t \upxt + \partial_x (\gamma \upxt) &= - \lambdaplus (\uplus, \uminus) \upxt + \lambdaminus (\uplus, \uminus) \umxt, \\
		\partial_t \umxt + \partial_x (- \gamma \umxt) &=  \lambdaplus (\uplus, \uminus) \upxt - \lambdaminus (\uplus, \uminus) \umxt, \\
		\upm(x, 0) &= \upm_0 (x), \quad x \in \R,
  \end{cases}
\end{align}
with the turning rates defined as
\begin{align*}
	\lambdapm (\uplus, \uminus) := \lambda_1 + \lambda_2 h (\ypm[\uplus, \uminus]).
\end{align*}
Here $\gamma$ is the constant speed, while the two constants $\lambda_1$ and $\lambda_2$ represent a baseline turning rate and a bias turning rate, respectively. For a biologically realistic case, the turning function $h$ should be a positive, increasing, and bounded functional that depends on the communication signals $\ypm$ perceived from neighbors. As in \cite{Eftimie-Vries-Lewis-2007} and \cite{Eftimie-Vries-Lewis-Lutscher-2007} we choose
\begin{align*}
	h (\ypm[\uplus, \uminus]) = 0.5 + 0.5 \tanh (\ypm[\uplus, \uminus] - y_0),
\end{align*}
where the constant $y_0$ is chosen such that for $\ypm[0] = 0$, the value of $\lambdapm(0)$ is determined only by $\lambda_1$. These signals $\ypm$ are emitted by neighbors moving to the right ($\uplus$) and to the left ($\uminus$):
\begin{align*}
	\begin{split}
		\ypm [\uplus, \uminus] := &q_r \int_{0}^{\infty} K_r(s) (u(x \pm s, t) - u(x \mp s, t)) ds\\
		- &q_a \int_{0}^{\infty} K_a(s) (u(x \pm s, t) - u(x \mp s, t)) ds\\
		+ &q_{al} \int_{0}^{\infty} K_{al}(s) (u^{\mp}(x \pm s, t) - \upm(x \mp s, t)) ds.
	\end{split}
\end{align*}
Here we define the total density as $\uxt = \upxt + \umxt$. The constants $q_r, q_a$, and $q_{al}$ represent the magnitudes of three social interactions: repulsion, attraction, and alignment, respectively. The interaction kernels $K_j$ (with $j = r, a, al,)$ are described by:
\begin{align*}
	K_j(s) := \dfrac{1}{\sqrt{2 \pi m_j^2}} \exp{\left(-\dfrac{(s -s_j)^2}{2 m_j^2}\right)}, \quad j = r, a, al,
\end{align*}
where $s_j, j = r, a, al,$ define the spatial regions for repulsive, alignment, and attractive
interactions, while $m_j := s_j/8$ define the width of these regions. We choose the
constants $m_j$ such that the support of more than $98\%$ of the mass of the kernels is inside the interval $[0, \infty)$. A more detailed description of this model can be found in \cite{Eftimie-Vries-Lewis-2007} and \cite{Eftimie-Vries-Lewis-Lutscher-2007}. \\

\textbf{Spatially homogeneous steady states.} For the rest of the paper (and in particular the numerical simulations), we assume that system \eqref{eqn:main} is defined on a bounded domain $[0, L]$ with periodic boundary conditions (to allow us to approximate the infinite domain by a finite domain). The non-local interaction terms are wrapped around the domain (see \cite{Eftimie-Vries-Lewis-2007} for further discussion). Since system \eqref{eqn:main} is conservative, let us define the total population density to be 
\begin{align*}
    A := \dfrac{1}{L} \int_0^L (\upxt + \umxt) dx.
\end{align*}
The spatially homogeneous steady states of \eqref{eqn:main} are the solutions $(\uplus, \uminus) = (u^*, A - u^*)$ of the steady-state equation
\begin{align} 
\begin{split}\label{eqn:steady-state}
    0 = &- u^* \left( \lambda_1 + 0.5 \lambda_2 + 0.5 \lambda_2 \tanh(A q_{al} - 2q_{al} u^* - y_0) \right) \\
    &+ (A - u^*) \left(\lambda_1 + 0.5 \lambda_2 + 0.5 \lambda_2 \tanh(-A q_{al} + 2q_{al} u^* - y_0) \right).
\end{split}
\end{align}
We note that $(\uplus, \uminus) = (A/2, A/2)$ is a solution to \eqref{eqn:steady-state} for all $\lambda_1, \lambda_2$ and $q_{al}$. In the case of $q_{al} = 0$, the steady-state equation \eqref{eqn:steady-state} has a unique solution $(\uplus, \uminus) = (A/2, A/2)$, and hence we have the unique spatially homogeneous steady state of \eqref{eqn:main}. \\

\begin{table}[!ht]
        \centering
	\caption{A list with the model parameters used during simulations.}
	\label{table:parameters}
	\begin{tabular}{|l|l|l|}
		\hline
		Parameter & Description & Fixed value \\ \hline
		$\gamma$ & Speed & $\gamma = 0.1$ \\ \hline
            $\lambda_1$ & Baseline turning rate & $\lambda_1 = 0.2$ \\ \hline
            $\lambda_2$ & Bias turning rate & $\lambda_2 = 0.9$ \\ \hline
            $y_0$ & Shift of the turning function & $y_0 = 2$ \\ \hline
            $q_a$ & Magnitude of attraction & $q_a = 1.1$ \\ \hline
            $q_r$ & Magnitude of repulsion & $q_r = 2.2$ \\ \hline
            $q_{al}$ & Magnitude of alignment & $q_{al} = 0.0$ \\ \hline
            $s_a$ & Attraction range & $s_a = 1$ \\ \hline
            $s_r$ & Repulsion range & $s_r = 0.25$ \\ \hline
            $s_{al}$ & Alignment range & $s_{al} = 0.5$ \\ \hline
            $m_a$ & Width of attraction kernel & $m_a = 1/8$ \\ \hline
            $m_r$ & Width of attraction kernel & $m_r = 0.25/8$ \\ \hline
            $m_{al}$ & Width of attraction kernel & $m_{al} = 0.5/8$ \\ \hline
            $A$ & Total population size & $A = 2$ \\ \hline
            $L$ & Size of bounded domain space $[0,L]$ & $L = 10$ \\ \hline
	\end{tabular}
\end{table}

\begin{figure}[!ht]
    \centering
    \resizebox{\textwidth}{!}{\input{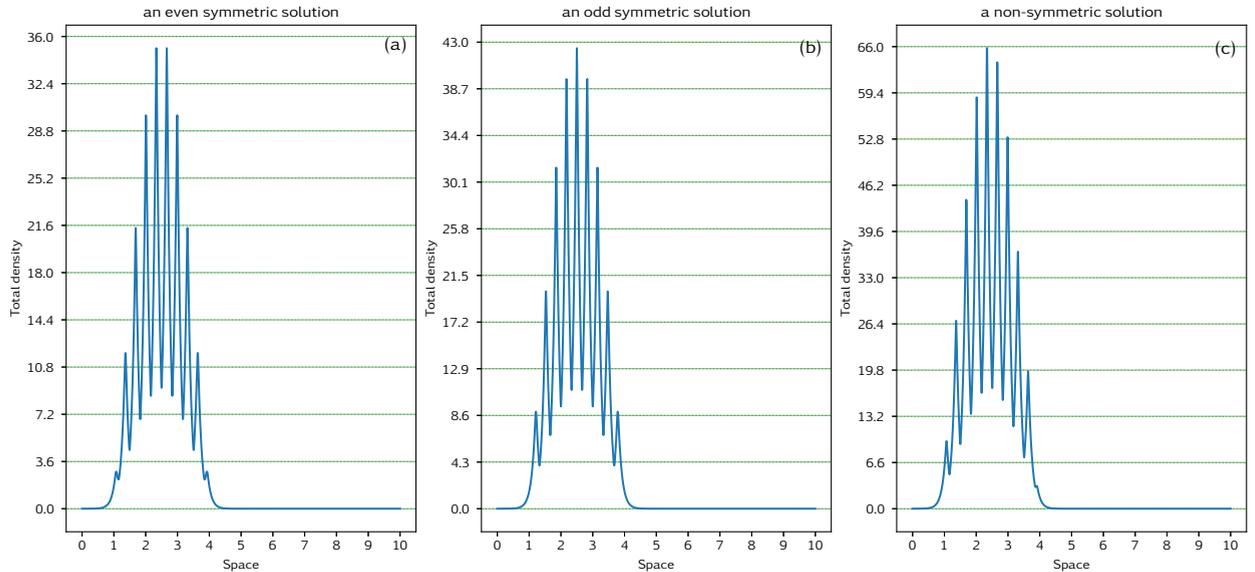}}
    \caption{Three types of solutions. These simulations used the upwind scheme and the initial condition \eqref{eqn:sin02_model}. \textbf{(a)} an even symmetric solution for $\hat{A}=3.5$. \textbf{(b)} an odd symmetric solution for $\hat{A}=5.0$. \textbf{(c)} a non-symmetric solution for $\hat{A}=10.0$. All other parameters are as in \cref{table:parameters}.}
    \label{fig:ODD_EVEN_NON}
\end{figure}

\textbf{Localized spatially-heterogeneous solutions with different symmetries.} 
In \cref{fig:Fi3_WeaklyNonlinear}, we observed that the stable branch with odd symmetric localized spatially-heterogeneous solutions can be obtained for $q_a \geq 1.095$, while the stable branch with even symmetric can be obtained for $q_a < 1.095$. Hence, it is possible to fix the value of $q_a = 1.1$, and consider the change of solutions symmetries by changing the initial amplitude. 
Indeed, for the parameter values listed in \cref{table:parameters} (and explained in more detail in~\cite[Table 1]{Eftimie-Vries-Lewis-Lutscher-2007} and \cite[Fig. 2 and Fig. 3]{Eftimie-Vries-Lewis-2009}), we can obtain (see~\cref{fig:ODD_EVEN_NON}) different localized solutions with different symmetries which depend on the amplitude of perturbations of the spatially homogeneous steady state $(\uplus, \uminus) = (1, 1)$.  For the initial condition 
\begin{align} \label{eqn:sin02_model}
    u(x, 0) = 2 \uplus(x, 0) = 2 \uminus(x, 0) = 2 + \hat{A} \, (0.5 + 0.5\sin(0.2 \pi x)), 
\end{align}
where $\hat{A}$ denotes the initial amplitude of the initial total density $u(x, 0)$, \cref{fig:ODD_EVEN_NON} illustrates three types of solutions: \textbf{(a)} for $\hat{A}=3.5$ we obtain an even symmetric solution consisting of 10 peaks; \textbf{(b)} for $\hat{A}=5.0$ we obtain an odd symmetric solution consisting of 9 peaks; \textbf{(c)} for $\hat{A}=10.0$ we obtain a non-symmetric solution. This suggests the possibility that these localized solutions belong to solution branches that are part of a snake-and-ladder bifurcation. To understand better these different types of localised solutions, in the next section, we consider different numerical schemes.

%% file: Numerical_scheme.tex
\section{Numerical schemes} \label{sec:numerical_scheme}

Given a finite time horizon $T > 0$, and a bounded domain space $[0, L]$, we consider the computation domain $[0, T] \times [0, L]$. Let $\Delta t$ and $\Delta x$ be the constant time and space steps, respectively. We set 
\begin{align*}
	N_t = \floor*{\dfrac{T}{\Delta t}}, \quad \text{and} \quad N_x = \floor*{\dfrac{L}{\Delta x}}. 
\end{align*}
Then for any $1 \leq i \leq N_x$, $0 \leq n \leq N_t$, we define the discrete mesh points $(x_i, t^n) = (i \Delta x, n \Delta t)$ and the cells $C_i = [x_{i - 1/2}, x_{i + 1/2})$. 
For $1 \leq i \leq N_x$, $1 \leq n \leq N_t$, we also denote by $(\uplus_{i})^n$, $(\uminus_{i})^n$ and $(u_{i})^n$ the approximation of the averages of $\uplus(x, t^n)$, $\uminus(x, t^n)$ and $u(x, t^n)$ on the cells $C_i$, namely
\begin{align*}
	(\uplus_{i})^n := \dfrac{1}{\Delta x} \int_{C_i} \uplus(x, t^n) dx, \quad (\uminus_{i})^n := \dfrac{1}{\Delta x} \int_{C_i} \uminus(x, t^n) dx, \quad u_{i}^n := \dfrac{1}{\Delta x} \int_{C_i} u(x, t^n) dx.
\end{align*}
As initial conditions we set, for $1 \leq i \leq N_x$,
\begin{align*}
	(\uplus_{i})^0 := \dfrac{1}{\Delta x} \int_{C_i} \uplus_0(x) dx, \qquad (\uminus_{i})^0 := \dfrac{1}{\Delta x} \int_{C_i} \uminus_0(x) dx. 
\end{align*}
At the boundary of the domain, we use periodic boundary conditions, namely
\begin{align*}
	\upm(x, t) = \upm(x - L) \quad \text{if } x > L, \quad \text{and} \quad \upm(x, t) = \upm(x + L) \quad \text{if } x < 0.
\end{align*}
This implies that for $0 \leq n \leq N_t$,
\begin{align*}
	(\upm_{i})^n = (\upm_{i - N_x})^n \quad \text{if } i > N_x, \quad \text{and} \quad (\upm_{i})^n = (\upm_{i + N_x})^n \quad \text{if } i < 0.
\end{align*}
We define the source terms in \eqref{eqn:main} as follows:
\begin{align*}
	(s_i^+)^n &:= - \lambda^+((\uplus_{i})^n, (\uminus_{i})^n) \, (\uplus_{i})^n + \lambda^-((\uplus_{i})^n, (\uminus_{i})^n) \, \uminus_{i})^n, \\ 
	(s_i^-)^n &:=  \lambda^+((\uplus_{i})^n, (\uminus_{i})^n) \, (\uplus_{i})^n - \lambda^-((\uplus_{i})^n, (\uminus_{i})^n) \, \uminus_{i})^n.
\end{align*}
To calculate these source terms, we approximate the infinite integrals by integrals on finite domains: $0 < s < 2s_j, \; j = r, a, al$; after that we approximate these integrals using the composite Simpson's 1/3 rule. This means that, 
\begin{align}
	\begin{split}
		Y_i^n &:= \int_{0}^{\infty} K_{al}(s) (\uminus(x_i + s, t_n) - \uplus(x_i - s, t_n)) ds\\
		&= \int_{0}^{2s_{al}} K_{al}(s) (\uminus(x_i + s, t_n) - \uplus(x_i - s, t_n)) ds \\
		&= \dfrac{\Delta x}{3} \left[K_{al}(0) \left[(\uminus_{i})^n - (\uplus_{i})^n\right] + K_{al}(2s_{al})\left[(\uminus_{i + N_{al}})^n - (\uplus_{i - N_{al}})^n\right] \right] \\
		&+ \dfrac{2 \Delta x}{3} \sum_{s = 1}^{N_{al}/2 -1} K_{al}(2s\Delta x)\left[(\uminus_{i + 2s})^n - (\uplus_{i - 2s})^n\right] \\
            &+ \dfrac{4 \Delta x}{3} \sum_{s = 1}^{N_{al}/2} K_{al}((2s-1)\Delta x)\left[(\uminus_{i + 2s-1})^n - (\uplus_{i - 2s+1})^n\right] 
	\end{split} \label{eqn:Y_i}
\end{align}
where $N_{al} := \frac{2s_{al}}{\Delta x}$. The same formula is applied for the remaining integrals. \\

Since the goal of this study is to investigate numerically the localized symmetric/asymmetric patterns that could suggest the presence of a snake-and-ladder bifurcation, in the following we discuss a few numerical schemes that have been used to obtain these localized solutions. To this end, we re-write \eqref{eqn:main} as follows:  
\begin{align} \label{eqn:main_numerical}
    \partial_t U + \partial_x F(U) = S(U),
\end{align}
where $U = [\uplus, \uminus]^{T}$, $F(U) = [\gamma \uplus, - \gamma \uminus]^{T}$, and $S(U) = [ - \lambdaplus  \uplus + \lambdaminus\uminus, \lambdaplus  \uplus - \lambdaminus \uminus ]^{T}$.

Throughout this study, we consider the following numerical schemes:
\begin{itemize}
\item \textbf{Upwind scheme.} We first apply the upwind scheme, which is equivalent to the Godunov method for our case where the velocity is constant ( see, e.g., \cite[Chapter 4]{LeVeque-2002-Finite-volume}). The discretized model reads
\begin{align*}
		(\uplus_{i})^{n+1} &= (\uplus_{i})^n - \dfrac{\gamma \Delta t}{\Delta x} \left[(\uplus_{i})^n - (\uplus_{i-1})^n\right] + \Delta t (s_i^+)^n, \\
		(\uminus_{i})^{n+1} &= (\uminus_{i})^n + \dfrac{\gamma \Delta t}{\Delta x} \left[(\uminus_{i+1})^n - (\uminus_{i})^n\right] + \Delta t (s_i^-)^n.
\end{align*}
\item \textbf{MacCormack scheme.} We consider a two-stage approach known as the MacCormack scheme (see \cite{MacCormack-2003}, see also \cite{Helbing-Treiber-1999}). The concept behind this scheme is to utilize upwind differences in the first stage and downwind differences in the second stage by using the values in the first stage. The values in the subsequent time step are calculated as the average of the previous step's values and the values from the second stage. The order in which the two directions are used can also be switched, or one can alternate between the two orderings in successive time steps, yielding a more symmetric method. Applying the MacCormack scheme, we obtain the first stage
\begin{align*}
		(\uplus_{i})^{*} &= (\uplus_{i})^n - \dfrac{\gamma \Delta t}{\Delta x} \left[(\uplus_{i})^n - (\uplus_{i-1})^n\right] + \Delta t (s_i^+)^n, \\
		(\uminus_{i})^{*} &= (\uminus_{i})^n + \dfrac{\gamma \Delta t}{\Delta x} \left[(\uminus_{i+1})^n - (\uminus_{i})^n\right] + \Delta t (s_i^-)^n,
\end{align*}
the second stage
\begin{align*}
		(\uplus_{i})^{**} &=  (\uplus_{i})^{*} - \dfrac{\gamma \Delta t}{\Delta x} \left[(\uplus_{i+1})^{*} - (\uplus_{i})^{*} \right] + \Delta t (s_i^+)^{*} , \\
		(\uminus_{i})^{**} &= (\uminus_{i})^{*} + \dfrac{\gamma \Delta t}{\Delta x} \left[(\uminus_{i})^{*} - (\uminus_{i-1})^{n+1/2}\right] + \Delta t (s_i^-)^{*} ,
\end{align*}
and finally:
\begin{align*}
		(\uplus_{i})^{n+1} &= \dfrac{1}{2} \left[ (\uplus_{i})^n + (\uplus_{i})^{**} \right], \\
		(\uminus_{i})^{n+1} &= \dfrac{1}{2} \left[ (\uminus_{i})^n + (\uminus_{i})^{**}\right].
\end{align*}

\item \textbf{Fractional-Step Method (FSM).} Next, we explore the Fractional-Step Method, see, e.g., \cite[Chapter 17]{LeVeque-2002-Finite-volume}. We first split the system \eqref{eqn:main_numerical} into two sub-problems: a homogeneous conservation law system and an ordinary differential system (ODEs) as follows
\begin{align*}
    &\text{Problem A:} \quad \partial_t U + \partial_x F(U) = 0, \\
    &\text{Problem B:} \quad \partial_t U = S(U).
\end{align*}
The idea behind the FSM is to discretize the original system by discretizing the two sub-problems in an alternating manner (by using standard methods for each sub-problem). 
Here, for the FSM approach, we use the upwind scheme for the homogeneous conservation law in problem A, and a two-stage Runge-Kutta method for ODEs (see, e.g., \cite[Chapter 5]{LeVeque-2007-Finite-difference}) in problem B. Hence, we obtain the A-step:
\begin{align*}
		(\uplus_{i})^{*} &= (\uplus_{i})^n - \dfrac{\gamma \Delta t}{\Delta x} \left[(\uplus_{i})^n - (\uplus_{i-1})^n\right], \\
		(\uminus_{i})^{*} &= (\uminus_{i})^n + \dfrac{\gamma \Delta t}{\Delta x} \left[(\uminus_{i+1})^n - (\uminus_{i})^n\right],
\end{align*}
the first stage of B-step:
\begin{align*}
		(\uplus_{i})^{**} &= (\uplus_{i})^{*} + \dfrac{\Delta t}{2} (s_i^+)^{*}, \\
		(\uminus_{i})^{**} &= (\uminus_{i})^{*} + \dfrac{\Delta t}{2} (s_i^-)^{*},
\end{align*}
and finally:
\begin{align*}
		(\uplus_{i})^{n+1} &= (\uplus_{i})^{*} + \Delta t \, (s_i^+)^{**}, \\
		(\uminus_{i})^{n+1} &= (\uminus_{i})^{*} + \Delta t \, (s_i^-)^{**}.
\end{align*}


\item \textbf{Quasi-Steady Wave-Propagation Algorithm (QSA).} We consider here the Quasi-Steady Wave-Propagation Algorithm (QSA) introduced in \cite{LeVeque-1998}. The basic idea of QSA is to introduce a new discontinuity in the center of each grid cell at the start of each time step, with a value $U_i^{L}$ on the left half of the cell and a value $U_i^{R}$ on the right half. These values are chosen so that
\begin{align} \label{eqn:average}
    \dfrac{1}{2} (U_i^{L} + U_i^{R}) = U_i,
\end{align}
and also, if possible, that 
\begin{align} \label{eqn:discrete_version}
    \dfrac{F(U_i^{R}) - F(U_i^{L})}{\Delta x} = S(U_i).
\end{align}
The condition \eqref{eqn:average} guarantees that the cell average is unchanged by the modiﬁcation, while \eqref{eqn:discrete_version}, if satisﬁed, means that the waves resulting from solving the Riemann problem at this new
discontinuity will exactly cancel the effect of the source term in this cell. Note that \eqref{eqn:discrete_version} is a discrete version of $\partial_x F(U) = S(U)$. So, it is not necessary to apply the source term any longer. One of the possible choices of $(U_i^{L}, U_i^{R})$ is that
\begin{align*}
    U_i^{L} = U_i - \delta_i \quad \text{and} \quad U_i^{R} = U_i + \delta_i,
\end{align*}
where $\delta_{i}$ satisfies \eqref{eqn:average} and \eqref{eqn:discrete_version} (for example, see \eqref{eqn:delta_i} for a detailed formula of $\delta_{i}$ for our problem \eqref{eqn:main}). 
Applying Godunov's method, we obtain that
\begin{align*}
    U_i^{n+1} = U_i^{n} - \dfrac{\Delta t}{\Delta x} \left[F(U_{i+1/2}^{*}) - F(U_{i-1/2}^{*}) \right],
\end{align*}
where the flux $F(U_{i+1/2}^{*})$ is determined by the modified values $(U_i^n)^{R}$ and $(U_{i+1}^n)^{L}$. Then, the high-resolution methods can be applied directly.

For our problem, by using QSA combining high-resolution methods, we have that
\begin{align*}
    (\uplus_{i})^{n+1} &= (\uplus_{i})^n - \dfrac{\gamma \Delta t}{\Delta x} \left[(\uplus_{i})^{nL} - (\uplus_{i-1})^{nR}\right] - \dfrac{1}{2} \dfrac{\gamma \Delta t}{\Delta x}(\Delta x - \gamma \Delta t) \left[(\sigma_{i}^+)^n - (\sigma_{i-1}^+)^n \right], \\
    (\uminus_{i})^{n+1} &= (\uminus_{i})^n + \dfrac{\gamma \Delta t}{\Delta x} \left[(\uminus_{i+1})^{nL} - (\uminus_{i})^{nR}\right] - \dfrac{1}{2} \dfrac{\gamma \Delta t}{\Delta x}(\Delta x - \gamma \Delta t) \left[(\sigma_{i+1}^-)^n - (\sigma_{i}^-)^n \right],
\end{align*}
where
\begin{align*}
    (\upm_{i})^{nL} = (\upm_i)^n - (\delta_i^{\pm})^n \quad \text{and} \quad (\upm_{i})^{nR} = (\upm_i)^n + (\delta_i^{\pm})^n
\end{align*}
with
\begin{align} \label{eqn:delta_i}
    (\delta_i^+)^n = \dfrac{\Delta x \, (s_i^+)^n}{2 \gamma} \quad \text{and} \quad  (\delta_i^-)^n = - \dfrac{\Delta x \, (s_i^-)^n}{2 \gamma}.
\end{align}
By using different slopes $(\sigma_{i}^{\pm})^n$ (for more information about slopes, please refer to \cite[Chapter 6]{LeVeque-2002-Finite-volume}), we obtain different schemes as detailed below. But first we denote
\begin{align}
    &\text{Centered slope:}  &(\sigma_{i}^{\pm})^n &= \dfrac{(u_{i+1}^{\pm})^n - (u_{i-1}^{\pm})^n - (\delta_{i+i}^{\pm})^n - 2 (\delta_{i}^{\pm})^n - (\delta_{i-i}^{\pm})^n}{2 \Delta x}, \label{eqn:centered_slope} \\
    &\text{Upwind slope:}  &(\sigma_{i}^{\pm})^n &= \dfrac{(u_{i}^{\pm})^n - (u_{i-1}^{\pm})^n - (\delta_{i}^{\pm})^n - (\delta_{i-i}^{\pm})^n}{\Delta x}, \label{eqn:upwind_slope} \\
    &\text{Downwind slope:}  &(\sigma_{i}^{\pm})^n &= \dfrac{(u_{i+1}^{\pm})^n - (u_{i}^{\pm})^n - (\delta_{i+1}^{\pm})^n - (\delta_{i}^{\pm})^n}{\Delta x}. \label{eqn:downwind_slope}
\end{align}
\begin{enumerate}[label=(\roman*)]
    \item \textbf{QSA scheme.} Here, we use no slope, i.e., $(\sigma_{i}^{\pm})^n = 0$. Hence, we get
    \begin{align*}
		(\uplus_{i})^{n+1} &= (\uplus_{i})^n - \dfrac{\gamma \Delta t}{\Delta x} \left[(\uplus_{i})^n - (\uplus_{i-1})^n\right] + \dfrac{\Delta t}{2} \left[(s_i^+)^n + (s_{i-1}^+)^n \right], \\
		(\uminus_{i})^{n+1} &= (\uminus_{i})^n + \dfrac{\gamma \Delta t}{\Delta x} \left[(\uminus_{i+1})^n - (\uminus_{i})^n\right] + \dfrac{\Delta t}{2} \left[ (s_{i+1}^-)^n + (s_{i}^-)^n \right].
    \end{align*}
    We observe that, for the linear system, QSA is equivalent to applying the upwind method to the original data but with a ﬁrst-order approximation to the source term also included. 
    \item \textbf{QSA\_Center scheme.} Here, we use the centered slope \eqref{eqn:centered_slope}.
    \item \textbf{QSA\_BW scheme.} Here, we follow the idea of the Beam-Warming method which gives a fully-upwind 3-point method. So, we use the upwind slope \eqref{eqn:upwind_slope} for $(\sigma_{i}^{+})^n$ and the downwind slope \eqref{eqn:downwind_slope} for $(\sigma_{i}^{-})^n$.
    \item \textbf{QSA\_LW scheme.} Here, we follow the idea of the Lax-Wendroff method which gives a centered 3-point method. So, we use the downwind slope \eqref{eqn:downwind_slope} for $(\sigma_{i}^{+})^n$ and the upwind slope \eqref{eqn:upwind_slope} for $(\sigma_{i}^{-})^n$.
    \item \textbf{QSA\_Minmod scheme.} Here, we use the \textit{minmod limiter} as follows
    \begin{align*}
        (\sigma_{i}^{\pm})^n = \text{minmod} \left[ \text{upwind slope } \eqref{eqn:upwind_slope}, \text{ downwind slope }\eqref{eqn:downwind_slope} \right],
    \end{align*}
    where
    \begin{align*}
        \text{minmod}(a,b) = 
        \begin{cases}
            \min(\abs{a}, \abs{b}) \quad &\text{if } ab > 0, \\
            0 &\text{if } ab \leq 0.
        \end{cases}
    \end{align*}
    \item \textbf{QSA\_Superbee scheme.} Here, we use the \textit{superbee limiter} as follows
    \begin{align*}
        (\sigma_{i}^{\pm})^n = \text{maxmod} \left( (\sigma_{i}^{\pm})^{(1)}, (\sigma_{i}^{\pm})^{(2)} \right),
    \end{align*}
    where
    \begin{align*}
        (\sigma_{i}^{\pm})^{(1)} &= \text{minmod} \left[ \text{upwind slope } \eqref{eqn:upwind_slope}, \text{ two times downwind slope }\eqref{eqn:downwind_slope} \right], \\
        (\sigma_{i}^{\pm})^{(2)} &= \text{minmod} \left[ \text{two times upwind slope } \eqref{eqn:upwind_slope},  \text{ downwind slope }\eqref{eqn:downwind_slope} \right],
    \end{align*}
    and 
    \begin{align*}
        \text{maxmod}(a,b) = 
        \begin{cases}
            \max(\abs{a}, \abs{b}) \quad &\text{if } ab > 0, \\
            0 &\text{if } ab \leq 0.
        \end{cases}
    \end{align*}
    \item \textbf{QSA\_MC scheme.} Here, we use the \textit{monotonized central-difference limiter} (MC limiter) as follows
    \begin{align*}
        (\sigma_{i}^{\pm})^n = \text{minmod} \left[ (\sigma_{i}^{\pm})^{(1)}, (\sigma_{i}^{\pm})^{(2)}, (\sigma_{i}^{\pm})^{(3)} \right],
    \end{align*}
    where
    \begin{align*}
        (\sigma_{i}^{\pm})^{(1)} &= \text{centered slope } \eqref{eqn:centered_slope},\\
        (\sigma_{i}^{\pm})^{(2)} &= \text{two times upwind slope } \eqref{eqn:upwind_slope}, \\
        (\sigma_{i}^{\pm})^{(3)} &= \text{two times downwind slope }\eqref{eqn:downwind_slope}.
    \end{align*}
\end{enumerate} 
    
\end{itemize}

\textbf{Order of accuracy of these schemes.} As known from any numerical analysis textbook, the upwind scheme is first-order accurate, while the MacCormack scheme is second-order accurate. The FSM scheme employed here is combining the first-order upwind scheme and the second-order Runge-Kutta method. In theory, the FSM scheme is first-order accurate. However, since the source terms are complicated, we expect the second-order ODEs solver to maintain overall accuracy. The QSA scheme is also first-order accurate. QSA\_Center, QSA\_BW, and QSA\_LW are second-order schemes. The high-resolution methods such as QSA\_Minmod, QSA\_Superbee, and QSA\_MC are formally not second-order accurate. However, the order of accuracy is not everything,  meaning that it is not always true that a method with a higher order of accuracy is more accurate on a particular grid or for a particular problem (for more information refer to \cite[Chapter 8]{LeVeque-2002-Finite-volume}).

\textbf{Calculating the error.} To investigate the convergence of the numerical solution to a steady state solution, we first define the discrete $\mathit{L^1}$ norm:
\begin{align*}
	\normLp{u(x, n \Delta t)}{1} = \normLp{u(x, t^n)}{1} := \dfrac{L}{N_x} \sum_{i = 1}^{N_x} \abs{u_i^n}.
\end{align*}
Then we denote by $E(t), \, t = 1,2,\ldots T,$ the discrete $\mathit{L^1}$-error between the total density $u$ at two adjacent time steps. To this end we take $t = k \Delta t (\in\{1,2,3...,T\})$ and calculate
\begin{align} \label{eqn:Et}
\begin{split}
    E(t) &:= \normLp{u(x,t)-u(x,t-\Delta t)}{1}= 
 \normLp{u(x, k \Delta t) - u(x, (k-1) \Delta t)}{1}\\ &= \normLp{u(x, t^k) - u(x, t^{k-1})}{1}
 = \dfrac{L}{N_x} \sum_{i = 1}^{N_x}  \abs{u_i^k - u_i^{k-1}}, \quad t = 1,2,\ldots T.
\end{split} 
\end{align}

\begin{remark} \label{remrak:L1-L2}
    \begin{enumerate}
        \item This paper employs the $\mathit{L^1}$-norm instead of the $\mathit{L^2}$-norm to quantify the numerical solution. There are several advantages associated with this choice. Firstly, the $\mathit{L^1}$-norm proves more suitable for solutions characterized by large distances between local maxima and minima. In particular, when the $\mathit{L^2}$-norm is utilized, a substantial discrepancy arises during the transition from odd symmetry to even symmetry. Furthermore, the $\mathit{L^2}$-norm of the solution can decrease when switching from even symmetry to odd symmetry, despite the increase in the $\mathit{L^2}$-norm of the initial data. Secondly, the $\mathit{L^1}$-norm of the total density remains independent of time within our model, facilitating control $\mathit{L^1}$-norm of solutions by adjusting the initial data. In fact, the $\mathit{L^1}$-norm of the total density can be approximated based on the initial amplitude $\hat{A}$ in \eqref{eqn:sin02_model} as 
        \begin{align} \label{eqn:L1-A}
            \normLp{u(x, t)}{1} \approx 20 + 5 \hat{A}, \qquad \forall t \geq 0.
        \end{align}
        \item Since we use periodic boundary conditions and need an approximation of the infinite integrals \eqref{eqn:Y_i} at each point on the grid, the adaptive mesh refinement technique seems challenging to apply. 
    \end{enumerate}
\end{remark}

%% file: Choice_numerical_scheme.tex
\section{Numerical aspects related to the identification of solutions with different symmetries} \label{sec:choice_scheme}
In this section, we discuss the advantages and disadvantages of the numerical schemes introduced in \cref{sec:numerical_scheme} that are relevant to the localised patterns observed in the non-local hyperbolic systems introduced in \cref{sec:model}. 
Since the localised solutions in Figure~\ref{fig:ODD_EVEN_NON} exhibit even/odd symmetry and asymmetry, we focus on the dependence of these different types of solutions on the initial condition, the space step, and the time step for each numerical scheme via various test cases. It should be noted that all simulations in this paper were executed using fixed parameters, as outlined in \cref{table:parameters}. 

We first state the stop criteria for our simulations (to ensure that the solution is approaching a steady state):
\begin{definition}[Stop criteria] A simulation will stop when it reaches
    \begin{enumerate}[label=(\roman*)]
    \item the given final time $T$, or
    \item the stop time $T^* := 1.34 \times t_0 < T$ where $t_0$ is the first time that satisfies $E(t_0) < \num{e-14}$, with the error $E(t)$ defined in \eqref{eqn:Et}.
    \end{enumerate}
\end{definition}
Next, we define two different types of numerical solutions that have a very small error $E(t)$ (which would suggest that we have reached a steady state, but this is not the case): 
\begin{definition}[Numerical transient solutions and numerical steady-state solutions] \label{def:sol}
A numerical solution at time $t^*$ characterized by a minimum error is called 
    \begin{itemize}
    \item  \textbf{transient solution} if $E(t^*) < \num{e-8}$ and $E(t)$ has a local minimum at $t^*$. 
    \item \textbf{steady-state solution} if $E(t^*) < \num{e-14}$ and $E(t)$ does not change significantly for $t > t^*$ (i.e., $E(t) < \num{e-14},$ $\forall t > t^*$).
    \end{itemize}
\end{definition}

\begin{remark}
	\begin{enumerate}
		\item Due to slow convergence to a steady-state solution, see \cref{subsec:transient_solution}, as well as the long execution time, see \cref{subsec:CPU_time}, this project requires the use of a fast-computation language. Therefore, we use \textbf{C++}, a low-level programming language with the \textbf{double} data type, which typically allows for 16 significant digits. In addition, we implemented the shared-memory parallel technique (OpenMP) to compute the integrals in the source terms (see formula \eqref{eqn:Y_i}).
        Besides, \textbf{Python}, a high-level programming language, is utilized for the post-processing of the results. 
		\item Numerical simulations are typically stopped at time $t^*$ if $E(t^*) < \num{e-8}$. However, for our specific problem, it is possible that the solution at time $t^*$ is only a transient solution and not the desired steady-state solution. In some cases, the value of $E(t)$ at a transient solution can be less than $\num{e-14}$ (see \cref{fig:transient_steady_state_solution}). Therefore, we set the stop time $T^* = 1.34 \, t_0$ where $t_0$ is the first time satisfied $E(t_0) < \num{e-14}$ in the stop criteria (ii) to observe the geometric of $E(t)$ for $t > t_0$. This allows us to observe the behavior of $E(t)$ for $t > t_0$, which serves as a basis for determining whether the numerical solution is a transient or steady-state solution. Since the limit of \textbf{double} data type in \textbf{C++}, it is not possible to use the condition $E(t_0) < \num{e-15}$ or less than.
		\item Typically, stop criteria (ii) are used to terminate iterative algorithms in our simulation, but in some cases, stop criteria (i) must be applied. These cases include situations where the value of the error $E(t)$ at a transient solution falls below $\num{e-14}$, as illustrated in \cref{fig:transient_steady_state_solution}, or when it is necessary to observe the behavior of solutions and the geometric of the error $E(t)$ over a prolonged period. In all other cases, stop criteria (ii) are sufficient for determining when to end the simulation.
	\end{enumerate}
\end{remark}

In \cref{subsec:transient_solution} we will discuss both transient and steady-state solutions, while \cref{subsec:CPU_time} will cover the execution time of simulations. In \cref{subsec:initial_data}, \cref{subsec:space_step}, and \cref{subsec:time_step} we will explore the impact of initial data, space step, and time step on the types of solutions, respectively. It's important to note that these subsections only address cases where the solution converges; non-convergence solutions will be discussed in \cref{subsec:non_convergence}. Finally, in \cref{subsec:overview}, we will provide the change in the solution symmetry as we vary the amplitude of initial conditions, for different numerical schemes.

{\bf Notation:} In the Tables that appear in the rest of this paper, we write, for short, {\textquoteleft ODD\textquoteright} for an odd symmetric solution, {\textquoteleft EVEN\textquoteright} for an even symmetric solution, and {\textquoteleft NON\textquoteright} for a non-symmetric solution.

\subsection{Numerical transient solutions vs. numerical steady-state solutions.} \label{subsec:transient_solution}

\begin{figure}[!ht] 
	\centering
        \resizebox{\textwidth}{!}{\input{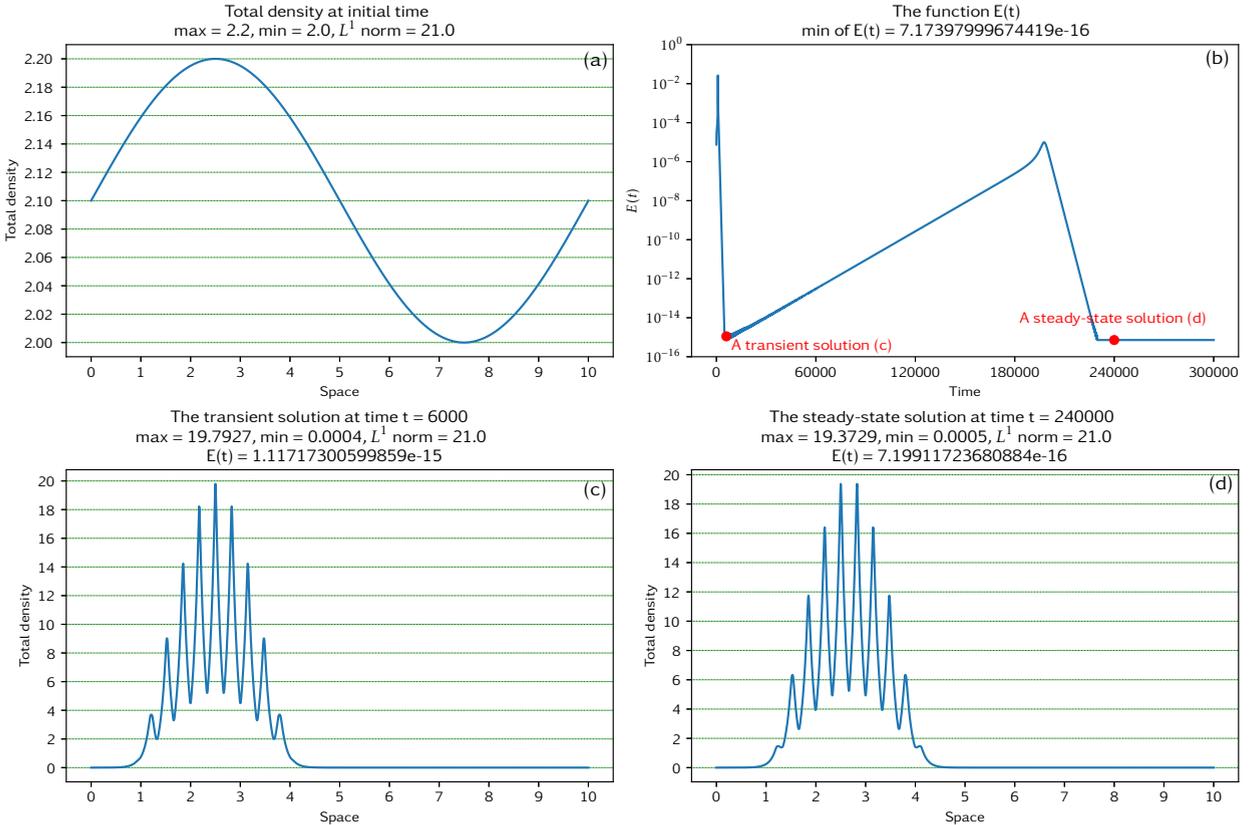}}
	\caption{A sinusoidal initial condition (panel (a)) could lead to an odd symmetric transient solution at $t = 6000$ (panel (c)) and an even symmetric steady-state solution at $T = 300000$ (panel (d)), both having minimal errors (panel (b))). The simulation used the upwind scheme with an initial amplitude of $\hat{A} = 0.2$, a space step of $\Delta x = 2^{-7}$ and a time step of $\Delta t = 2^{-6}$. }
	\label{fig:transient_steady_state_solution}
\end{figure}

One of the numerical issues encountered so far is the extremely slow convergence of the transient solution to a steady-state solution that can be followed numerically through the continuation algorithm, see \cite{Uecker-2022}. We discuss this in more detail for our particular problem, the nonlocal hyperbolic system \eqref{eqn:main}. 

We begin by showing in \cref{fig:transient_steady_state_solution} a simulation obtained using the upwind scheme, where the amplitude of the initial condition is $\hat{A} = 0.2$ (see sub-panel (a)), the spatial step is $\Delta x = 2^{-7}$, and the time step is $\Delta t = 2^{-6}$. The error $E(t)$ calculated up to time 300000 (see sub-panel (b)) can achieve minimal values at different time points: at an early time $t=6000$ (see the transient solution in \cref{fig:transient_steady_state_solution}(c)), and at some later times $t>240000$ (see the steady-state solution in \cref{fig:transient_steady_state_solution}(d)). 
However, the type of transient solution (odd symmetric) differs from that of a steady-state solution (even symmetric).
Note that these two numerical solutions are likely the result of the fact that the spatially heterogeneous steady state (i.e., the localised solution) is an unstable saddle point. The numerical solution trajectories approach quickly the transient solution along the stable manifold (due to a large negative eigenvalue). Then, these solutions move very slowly away from this transient state along the unstable manifold of the saddle point (due to a very small positive eigenvalue). For a very large time $t>200000$ in \cref{fig:transient_steady_state_solution} we see that the solution eventually approaches a stable spatially heterogenous state (steady-state solution). Besides, the total density at the local maximum of $E(t)$ is non-symmetric, see \cref{fig:transient_at_some_time} in Appendix A.1.

As shown in \cref{fig:transient_steady_state_solution}, we have a transition from odd symmetric to even symmetric solutions; we also observe the opposite transition: from even symmetric to odd symmetric solutions in some other cases. Additionally, in certain cases, we could also observe a transition from either odd or even symmetric to non-symmetric and the opposite transition. The transient solutions were observed to exist for all schemes analyzed in this paper, with varying initial amplitudes and different time and space steps. For further information, see Appendix A.1.

The presence of the transient solution presents us with two challenges. First, the small error $E(t)$ associated with the transient solutions can cause confusion and misinterpretation of these solutions as steady-state solutions. Hence, using the transient solution results in a distinct branching pattern (with odd symmetry) compared to the branching pattern (with even symmetry) obtained with the steady-state solution. Second, the convergence rate in our simulations is considerably slow, necessitating a significantly higher final time value to attain the steady-state solutions. Notably, the convergence rate is slower when the error of the transient solutions is smaller.

It was expected that decreasing the space or time step would eliminate transient solutions or at least improve the convergence rate by increasing the value of $E(t)$ at the transient solutions. However, this was not always the case. In some cases, reducing the space step could lead to the creation of a transient solution (which will be discussed in more detail in \cref{subsec:space_step}), while decreasing the time step could reduce the value of $E(t)$ at the transient solution, thereby slowing down the convergence (which will be discussed in more detail in \cref{subsec:time_step}). 

\subsection{Execution time of different numerical schemes.} \label{subsec:CPU_time} 
\begin{figure}[!ht] 
	\centering
        \begin{adjustbox}{clip,trim=0.8cm 0.7cm 0cm 0.8cm}
        \resizebox{1.1\textwidth}{!}{\input{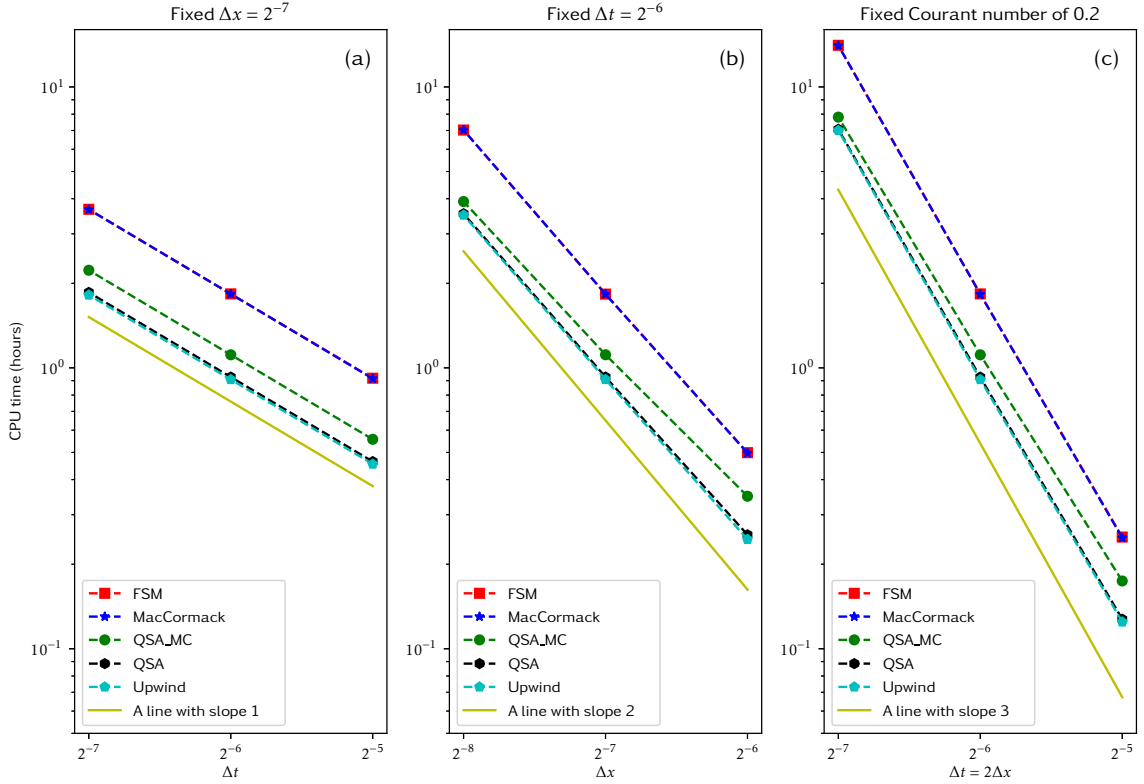}}
        \end{adjustbox}
	\caption{The execution time (CPU time) for simulations with the same final time $T = 1000$. Computations have been performed on the supercomputer facilities of the Mésocentre de Calcul at the University of Franche-Comté.}
	\label{fig:CPU_time}
\end{figure}

To compare the results obtained with different numerical schemes, we discuss the execution (CPU) time for our simulations with a fixed final time of $T = 1000$. 
In \cref{fig:CPU_time} we first compare the CPU time for a space step of $\Delta x = 2^{-7}$ and a time step of $\Delta t = 2^{-6}$: the upwind and QSA schemes are the fastest, with a CPU time of around 0.91 hours; the QSA with slope schemes is slower, taking around 1.11 hours; the two-stage schemes, MacCormack and FSM, are the slowest, with a CPU time of around 1.83 hours, see \cref{fig:CPU_time}(a).

If we fix the value of the space step $\Delta x$, then decreasing the time step $\Delta t$ by a factor of two will increase the CPU time by a factor of two, see \cref{fig:CPU_time}(a). On the other hand, if we fix the value of the time step $\Delta t$, decreasing the space step $\Delta x$ by a factor of two will increase the CPU time by a factor of four, see \cref{fig:CPU_time}(b).
Therefore, if we fix the Courant number at 0.2, i.e., $\Delta t = 2 \Delta x$, then increasing both the space and time steps by a factor of two will increase the CPU time by a factor of eight, see \cref{fig:CPU_time}(c).

Due to the slow convergence towards steady-state solutions, the final time or stop time required is often large, usually ranging from $4000$ to $16000$ in cases with no transient solution and a space step of $\Delta x = 2^{-7}$. If the space step is smaller, the final time may be even larger, as shown in \cref{fig:test_dx_error}. In cases where a transient solution exists, the final time required is very large, typically exceeding $50000$. For instance, in \cref{fig:transient_steady_state_solution}, a final time of $t > 200000$ was necessary to obtain the steady-state solution. Consequently, the CPU time required for our simulations is also very large, it takes 206 hours. Another example, when using the MacCormack scheme with a space step of $\Delta x = 2^{-9}$ and a time step of $\Delta t = 2^{-7}$, it took $2743$ CPU hours to converge to a stop time of $T^* = 49292$. Simulations with the space step of $\Delta x = 2^{-10}$ cannot be executed due to insufficient computational resources (it passes the limit of computation time of eight days for 16 cores on the supercomputer facilities of the Mésocentre de Calcul at the University of Franche-Comté). Since the identification of localised solutions with specific symmetries (that could form snake-and-ladder bifurcation branches) necessitates a significant number of simulations, the selection of numerical schemes as well as their corresponding parameters such as the time and space steps, or the selection of initial conditions, become crucial for the simulations as will be discussed in the following subsections.

\subsection{The influence of the initial data. }\label{subsec:initial_data} 

\begin{figure}[!ht]
	\centering
	\resizebox{\textwidth}{!}{\input{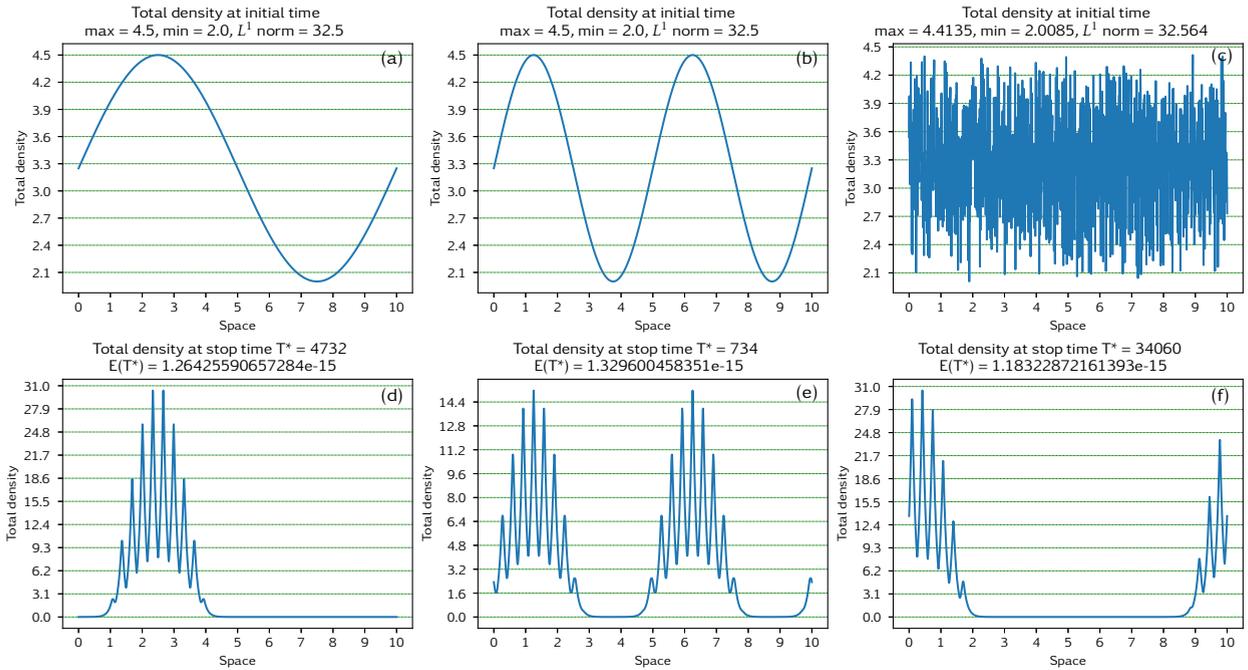}}
        \caption{The influence of the initial data on the types of solutions. The first column presents an even symmetric solution using the initial condition \eqref{eqn:sin02}. The second column presents a solution that includes two odd symmetric aggregations using the initial condition \eqref{eqn:sin04}. The last column presents a non-symmetric solution using the initial condition \eqref{eqn:rand}. Three above solutions were obtained by using an initial amplitude of $\hat{A} = 2.5$, a space step of $\Delta x = 2^{-7}$, a time step of $\Delta t = 2^{-6}$ and the upwind scheme.}
	\label{fig:different_initial_conditions}
\end{figure}

Throughout this study we consider initial conditions that are perturbations of the spatially homogeneous steady state $(\uplus, \uminus) = (1, 1)$. In the following we discuss the following three initial conditions:
\begin{align} \label{eqn:sin02}
	u(x, 0) = 2 \uplus(x, 0) = 2 \uminus(x, 0) = 2 + \hat{A} \, (0.5 + 0.5\sin(0.2 \pi x)),
\end{align}
\begin{align}\label{eqn:sin04}
	u(x, 0) = 2 \uplus(x, 0) = 2 \uminus(x, 0) = 2 + \hat{A} \, (0.5 + 0.5\sin(0.4 \pi x)),
\end{align}
\begin{align}\label{eqn:rand}
	u(x, 0) = 2 \uplus(x, 0) = 2 \uminus(x, 0) = 2 + \hat{A} \, rand([0,1)),
\end{align}
where $rand([0,1))$ is a random number in $[0,1)$ and $\hat{A}$ denotes the amplitude of the initial perturbation. \cref{fig:different_initial_conditions}(a), (b), and (c) shows these initial conditions for an initial perturbation of amplitude $\hat{A} = 2.5$ (as well as a space step $\Delta x = 2^{-7}$ and a time step $\Delta t = 2^{-6}$). Note that the shape of solutions generated by initial data in the form of $\sin(\cdot)$ or $\cos(\cdot)$ is identical; here we choose the $\sin(\cdot)$ form, as the aggregations of peaks of the solution lie in the middle of the domain, making them easier to visualize. 


\begin{table}[!ht]
        \centering
	\caption{The influence of the initial data on the types of solutions. These simulations were executed with an initial amplitude of $\hat{A} = 2.5$, a space step of $\Delta x = 2^{-7}$, and a time step of $\Delta t = 2^{-6}$. We denote 2-ODD as a solution that includes two odd symmetric aggregations.}
	\label{table:different_initial_conditions}
	\begin{tabular}{|l|c|c|c|}
		\hline
		\diagbox[width=12em]{Scheme used}{Initial condition} & \eqref{eqn:sin02} & \eqref{eqn:sin04} & \eqref{eqn:rand} \\ \hline
		Upwind & EVEN & 2-ODD & NON 	\\ \hline
		MacCormack & EVEN & 2-ODD & NON	\\ \hline
		FSM & EVEN & 2-ODD & NON 	\\ \hline
		QSA & EVEN & 2-ODD & NON 	\\ \hline
		QSA\_Center & EVEN & 2-ODD & NON 	\\ \hline
		QSA\_BW & EVEN & 2-ODD & NON 	\\ \hline
		QSA\_LW & EVEN & 2-ODD & NON\\ \hline
		QSA\_Minmod & ODD & 2-ODD & NON 	\\ \hline
		QSA\_Superbee & EVEN & 2-ODD & EVEN	\\ \hline
		QSA\_MC & ODD & 2-ODD & NON 	\\ \hline
	\end{tabular}
\end{table}

In \cref{fig:different_initial_conditions}(d) we observe that using the initial condition \eqref{eqn:sin02} (depicted in \cref{fig:different_initial_conditions}(a)), we obtain an even symmetric solution. However, using the initial condition \eqref{eqn:sin04} (as in \cref{fig:different_initial_conditions}(b)), we obtain a solution that includes two odd symmetric aggregations, see \cref{fig:different_initial_conditions}(e). \cref{fig:different_initial_conditions}(f) shows a non-symmetric solution where the initial condition \eqref{eqn:rand} is used. 
For a summary of the results obtained with various numerical schemes, for the initial conditions \eqref{eqn:sin02}, \eqref{eqn:sin04}, and \eqref{eqn:rand}, with an initial amplitude $\hat{A} = 2.5$, a space step $\Delta x = 2^{-7}$, and a time step $\Delta t = 2^{-6}$, see \cref{table:different_initial_conditions}. For the initial condition \eqref{eqn:sin02}, we obtain an even symmetric solution, except for the QSA\_Minmod and QSA\_MC schemes, which result in an odd symmetric solution. On the other hand, using the initial condition \eqref{eqn:sin04}, we obtain a solution that includes two odd symmetric aggregations. Finally, we obtain a non-symmetric solution using the initial condition \eqref{eqn:rand}, except for the QSA\_Superbee scheme which results in an even symmetric solution. 

It is worth noting that in the random case, obtaining a steady-state solution takes considerably longer, leading to slow convergence compared to other cases. Therefore, in the rest of this paper, we will use the initial condition \eqref{eqn:sin02} for all simulations.

\subsection{The influence of the space step. } \label{subsec:space_step}

\begin{figure}[!ht] 
	\centering
        \begin{adjustbox}{clip,trim=1.5cm 0.9cm 0cm 0.8cm}
        \resizebox{1.2\textwidth}{!}{\input{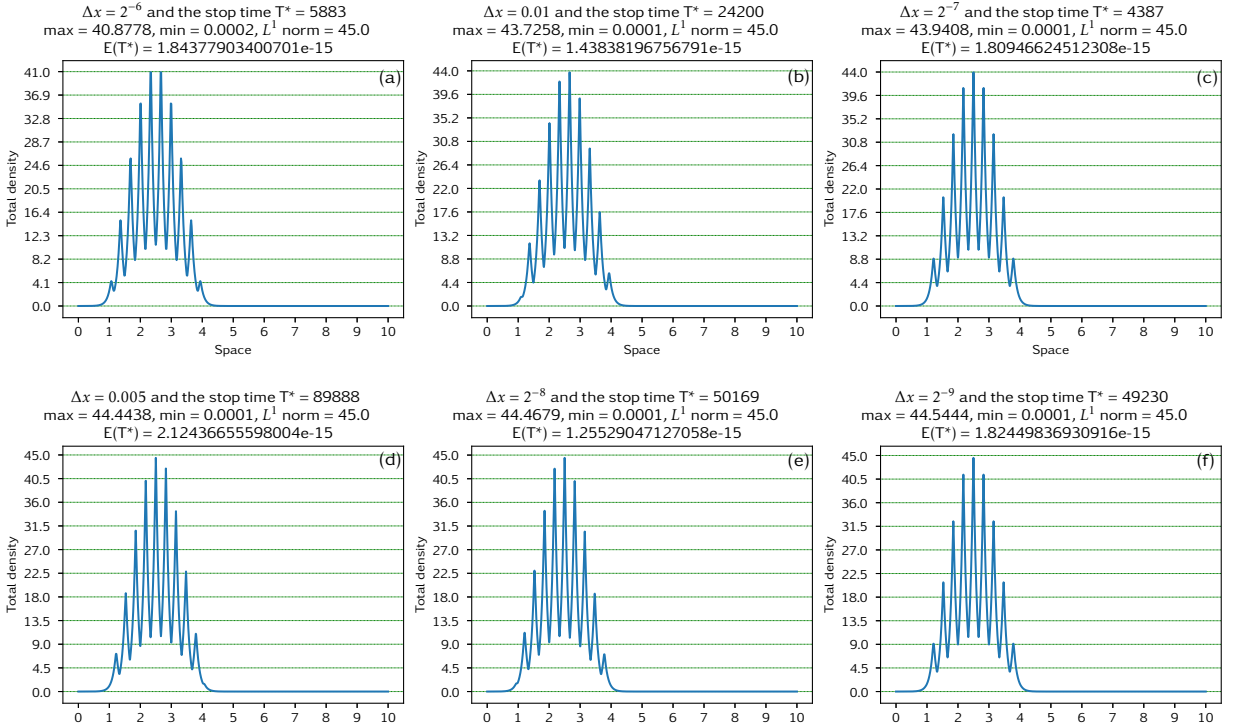}}
        \end{adjustbox}
	\caption{The influence of the space step on the types of solutions. These simulations were executed by using the QSA\_MC scheme and the initial condition \eqref{eqn:sin02} with an initial amplitude of $\hat{A} = 5.0$, and a time step of $\Delta t = 2 \Delta x$. In this case, it is possible to obtain all three types of solutions by using different space steps: an even symmetric solution in (a); odd symmetric solutions in (c) and (f); non-symmetric solutions in (b), (d) and (e).}
	\label{fig:test_dx_MC}
\end{figure}

Now, we discuss the influence of the space step on the types of solutions. To this end, we perform simulations with different space steps $\Delta x = \{2^{-6}; \, 0.01; \, 2^{-7}; \, 0.005; \, 2^{-8}; \, 2^{-9}\}$, with initial condition \eqref{eqn:sin02} having a fixed initial amplitude of $\hat{A} = 5.0$, and a fixed the Courant number of 0.2, i.e, $\Delta t = 2 \Delta x$. 
Numerical results show that using different space steps (while keeping all model parameters fixed) could lead to localised solutions with different symmetries. For example, \cref{fig:test_dx_MC} presents three types of solutions obtained with the QSA\_MC scheme: even symmetric solutions when $\Delta x = 2^{-6}$, see sub-panel (a); odd symmetric solutions when $\Delta x = 2^{-7}$ and $\Delta x = 2^{-9}$ see sub-panels (c) and (f); non-symmetric solutions when $\Delta x = 0.01$, $\Delta x = 0.005$ and $\Delta x = 2^{-8}$ see sub-panels (b), (d) and (e).
We also obtain all three types of solutions using QSA\_Minmod. For other schemes, it is possible to obtain two types of solutions, which are odd symmetric solutions and non-symmetric solutions. For more details, please refer to \cref{table:different_space_step}.
From the \cref{table:different_space_step}, it is also observed that different solution symmetries are obtained by different schemes, despite using the same space step in the simulations. For instance, all odd symmetric solutions are obtained for the space step $\Delta x = 2^{-9}$, while all non-symmetric solutions are obtained for the space step $\Delta x = 2^{-8}$. However, for the space steps $\Delta x = 0.005$, $\Delta x = 2^{-7}$, and $\Delta x = 0.01$, both odd symmetric and non-symmetric solutions can be obtained. The space step $\Delta x = 2^{-6}$ yields all three types of solutions.

\begin{table}[!ht]
        \centering
	\caption{The influence of the space step on the types of solutions. These simulations were executed by using the initial condition \eqref{eqn:sin02} with an initial amplitude of $\hat{A} = 5.0$, and a time step of $\Delta t = 2 \Delta x$.}
	\label{table:different_space_step}
	\begin{tabular}{|l|c|c|c|c|c|c|}
		\hline
		\diagbox[width=12em]{Scheme used}{$\Delta x$} & $2^{-6}$ & $0.01$ & $2^{-7}$ & $0.005$ & $2^{-8}$ & $2^{-9}$	\\ \hline
		Upwind & ODD & NON  & ODD & NON & NON & ODD	\\ \hline
		MacCormack & ODD & ODD & ODD & NON & NON & ODD \\ \hline
		FSM & ODD & NON  & ODD & NON & NON & ODD	\\ \hline
		QSA & NON & ODD & ODD & NON & NON & ODD \\ \hline
		QSA\_Center & NON & ODD & ODD & NON & NON & ODD \\ \hline
		QSA\_BW & NON & ODD & ODD & NON & NON & ODD \\ \hline
		QSA\_LW & NON & ODD & ODD & NON & NON & ODD \\ \hline
		QSA\_Minmod & EVEN & NON & NON & NON & NON & ODD \\ \hline
		QSA\_Superbee & NON & ODD & ODD & ODD & NON & ODD \\ \hline
		QSA\_MC & EVEN & NON & ODD & NON & NON & ODD \\ \hline
	\end{tabular}
\end{table}

\begin{figure}[!ht] 
	\centering
        \begin{adjustbox}{clip,trim=1.3cm 0.4cm 0cm 0.cm}
        \resizebox{1.2\textwidth}{!}{\input{Figures/pp1_test_dx_error.pgf}}
        \end{adjustbox}
	\caption{The influence of the space step. The function E(t) corresponds to varying space steps for several schemes: \textbf{(a)} the upwind scheme, \textbf{(b)} the MacCormack scheme, \textbf{(c)} the QSA\_BW scheme, \textbf{(d)} the QSA\_Minmod scheme, \textbf{(e)} the QSA\_Superbee scheme and \textbf{(f)} the QSA\_MC scheme. These simulations were executed by using the initial condition \eqref{eqn:sin02} with an initial amplitude of $\hat{A} = 5.0$, and a time step of $\Delta t = 2 \Delta x$.} 
	\label{fig:test_dx_error}
\end{figure}

In regard to the impact of the decreasing the space step on the convergence rate of the steady-state solutions, we observe that the smaller space step does not always converge faster than the large space step. In some cases, reducing the space step could lead to the creation of a transient solution.
For example, \cref{fig:test_dx_error} shows the error function $E(t)$ for several numerical schemes: the upwind (sub-panel (a); same properties for the FSM scheme), the MacCormack (sub-panel (b)), the QSA\_BW (sub-panel (c); same properties for the QSA, QSA\_Center, and QSA\_LW schemes), the QSA\_Minmod (sub-panel (d)), the QSA\_Superbee (sub-panel (e)), and the QSA\_MC (sub-panel (f)).  In these simulations, it is particularly noticeable that smaller space steps, such as $\Delta x = {0.005; \, 2^{-8}; \, 2^{-9}}$, consistently exhibit slower convergence compared to larger space steps like $\Delta x = {2^{-6}; \, 0.01; \, 2^{-7}}$. Additionally, it is observed that transient solutions do not arise for large space steps, such as $\Delta x = {2^{-6}; \, 0.01; \, 2^{-7}}$, whereas they do occur for smaller space steps.

In accordance with the theory of convergent numerical methods (for example, see \cite[Chapter 2]{LeVeque-2007-Finite-difference}), decreasing the space step (while maintaining the Courant number) should lead to numerical solutions that approach the exact solution. The observed changes in solution types, as described earlier, raise the question of convergence in the theoretical framework of the investigated numerical schemes for the non-local hyperbolic system. Furthermore, due to the slow convergence rate, selecting an appropriate space step becomes a challenging task in simulations. To strike a balance between achieving a satisfactory convergence rate and maintaining an acceptable computation time, as discussed in Section \ref{subsec:CPU_time}, we opted for a space step of $\Delta x = 2^{-7}$. This choice corresponds to a total of 1281 points on the spatial grid for all subsequent simulations presented in this paper.

\begin{remark}
    Note that it is difficult to simulate with smaller space steps, such as $2^{-10} = 0.000976563$, as the execution time would be excessively  (at least on the servers of supercomputer facilities of the Mésocentre de Calcul at the University of Franche-Comté), see \cref{subsec:CPU_time}. It is also difficult to use the adaptive mesh refinement technique, due to the changes in the solution symmetries if we change the space step and necessary approximations of the infinite integrals, see \cref{remrak:L1-L2}. Therefore, we conclude that such choices regarding the numerical schemes (and corresponding numerical space-time steps) depend also on the local computing environment at different academic institutions.
\end{remark}

\subsection{The influence of the time step. } \label{subsec:time_step} 
\begin{figure}[!ht] 
	\centering
        \resizebox{0.82\textwidth}{!}{\input{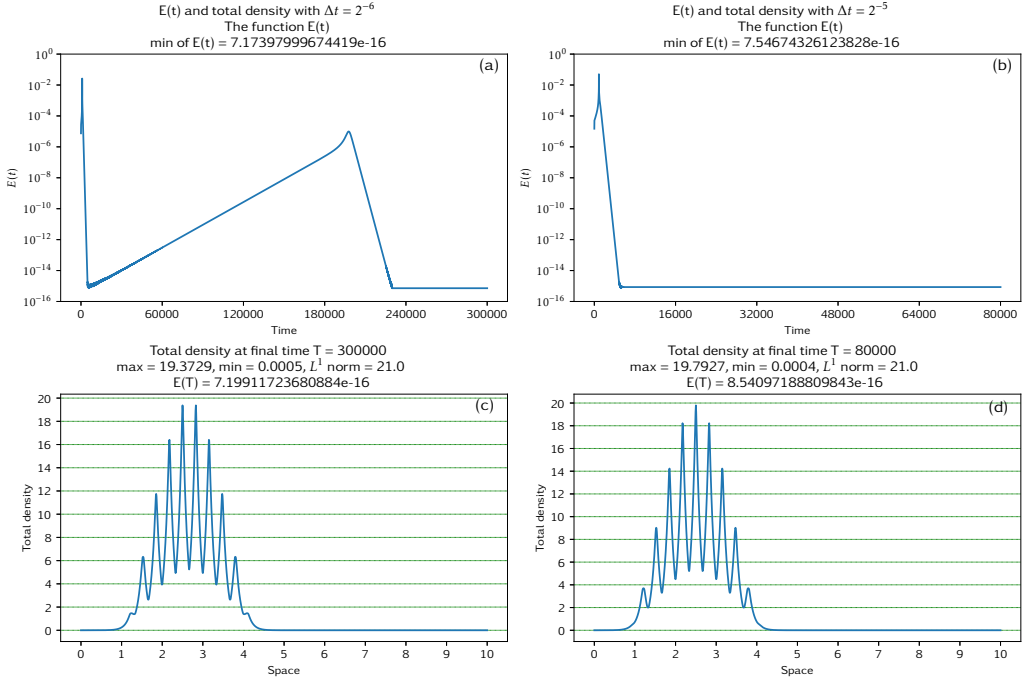}}
	\caption{The influence of the time step on the types of solutions. In the first column, an even symmetric solution is obtained with a time step of $\Delta t = 2^{-6}$, however, in the second column, an odd symmetric solution is obtained by using a time step of $\Delta t = 2^{-5}$. There exists a transient solution in the case of the time step of $\Delta t = 2^{-6}$. These simulations were executed by using the upwind scheme, and the initial condition \eqref{eqn:sin02} with an initial amplitude of $\hat{A} = 0.2$, and a space step of $\Delta x = 2^{-7}$.}
	\label{fig:test_dt_0200_upwind}
\end{figure}

Since we are seeking localised steady-state solutions, we expect that the solutions will remain constant across different time steps, and thus the $\mathit{L}^1$-norm of the difference between two solutions corresponding to consecutive time steps should be small. In this subsection, we discuss the influence of the time steps on such steady-state solutions, when we use the initial condition \eqref{eqn:sin02} and a fixed space step of $\Delta x = 2^{-7}$. 

In \cref{fig:test_dt_0200_upwind} we examine two simulations obtained with the upwind scheme, for an initial amplitude of $\hat{A} = 0.2$. The first simulation is obtained with a time step of $\Delta t = 2^{-6}$ and yields an even symmetric solution (\cref{fig:test_dt_0200_upwind}(c)). The second simulation, which employs a time step of $\Delta t = 2^{-5}$, produces an odd symmetric solution (\cref{fig:test_dt_0200_upwind}(d)). Note that, for the time step $\Delta t = 2^{-6}$, there exists also a transient solution; see \cref{fig:test_dt_0200_upwind}(a). However, it is important to note that differences in the type of solutions generated by varying time steps also exist in cases where there are no transient solutions; see Appendix A.2. For a summary of the localised patterns with different symmetries obtained with two different time steps $\Delta t = 2^{-5}$ and $\Delta t = 2^{-6}$, across all numerical schemes (using 361 different initial amplitude values ranging from $0.001$ to $36$), please see \cref{fig:all_schemes_7_5_6}. There we observe that variations in the symmetry of solutions due to changes in the time step occurs for all the schemes examined in this study, albeit at different values of initial amplitude $\hat{A}$. 
For further details, please see \cref{subsec:overview}. Also, for a further discussion on the impact of different time steps on model solution, see Appendix A.2.



\begin{figure}[!ht] 
	\centering
        \begin{adjustbox}{clip,trim=1.0cm 0.6cm 0cm 0.cm}
        \resizebox{\textwidth}{!}{\input{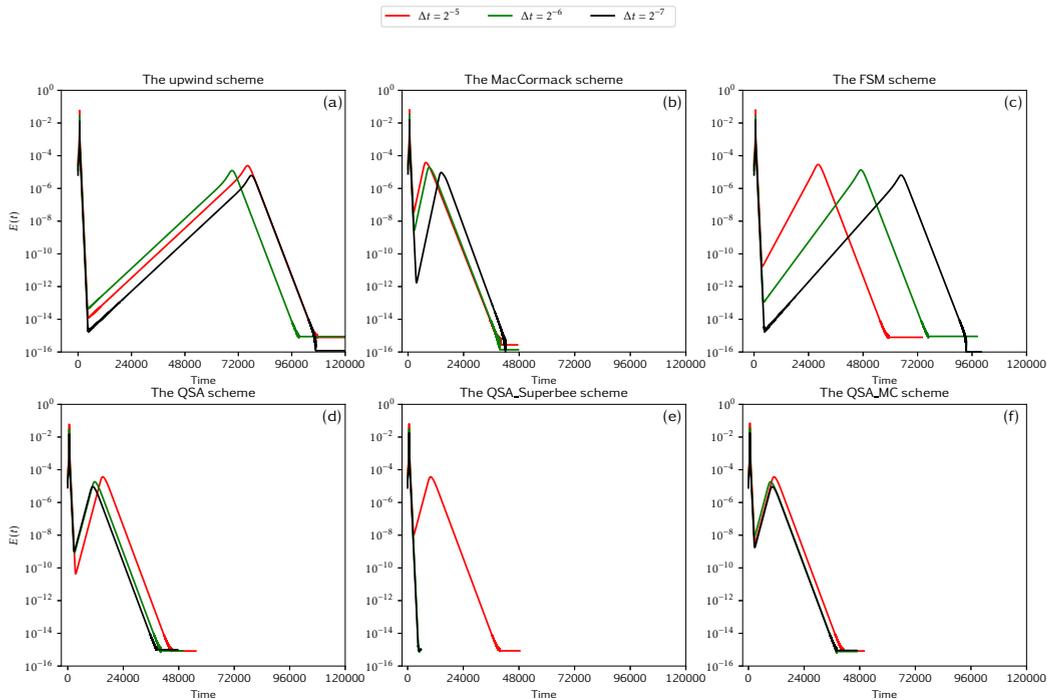}}
        \end{adjustbox}
	\caption{The error function $E(t)$ for varying time steps $\Delta t = \{2^{-5}; \, 2^{-6}; \, 2^{-7} \}$. The simulation considered an initial condition given by \eqref{eqn:sin02}, with the amplitude of perturbations $\hat{A} = 0.3$ and a space step $\Delta x = 2^{-7}$. }
	\label{fig:0300_some_schemes}
\end{figure}

Now we discuss the impact of decreasing the time step on the convergence rate of the steady-state solutions. To elaborate further, \cref{fig:0300_some_schemes} illustrates the behavior of $E(t)$ with varying time steps $\Delta t = 2^{-5}, \Delta t = 2^{-6}, \Delta t = 2^{-7}$, using an initial amplitude of $\hat{A} = 0.3$ and a space step $\Delta x = 2^{-7}$. This figure reveals that the MacCormack and the FSM schemes converge more slowly as the time steps decrease. These are likely a result of the fact that these schemes are two-stage schemes.  
Particularly for the upwind and FSM schemes with a time step of $\Delta t = 2^{-7}$, the values of $E(t)$ at transient solutions are less than $\num{e-14}$, and it takes over $t > 100000$ to reach steady-state solutions. However, in the case of the QSA\_Superbee scheme, the transient solution only occurs for the time step $\Delta t = 2^{-5}$, making it the most efficient scheme in this regard. For the other schemes, the convergence rate is slightly faster with smaller time steps, but transient solutions still persist. 

As in the previous subsection, the results here also raise the question of convergence in the theoretical framework of the investigated numerical schemes for nonlocal hyperbolic systems, and the difficulty to choose an appropriate time step for the simulations. 

\subsection{The non-convergence of numerical schemes. }\label{subsec:non_convergence}
\begin{figure}[!ht] 
	\centering
        \resizebox{\textwidth}{!}{\input{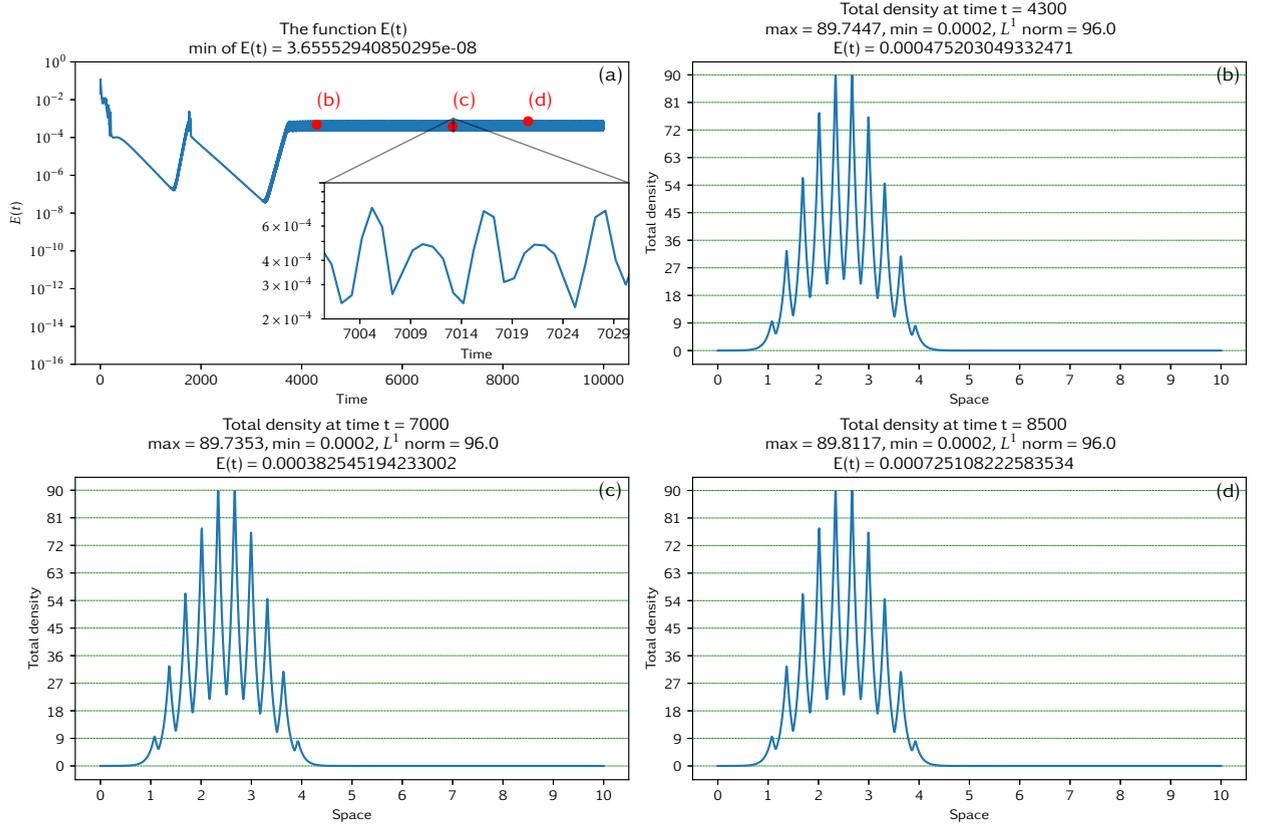}}
	\caption{Non-convergence of the FSM scheme using the initial condition \eqref{eqn:sin02} with an initial amplitude of $\hat{A} = 15.2$, a space step of $\Delta x = 2^{-7}$, and a time step of $\Delta t = 2^{-7}$. \textbf{(a)} the error function $E(t)$ and a subfigure is a zoom-in for an interval of time from $t = 7000$ to $t = 7300$; \textbf{(b)} a non-symmetric total density at time $t = 4301$ with $E(t) = 0.000475$; \textbf{(c)} a non-symmetric total density at time $t = 7001$ with $E(t) = 0.000382$. \textbf{(d)} a non-symmetric total density at time $t = 8501$ with $E(t) = 0.000725$.}
	\label{fig:non_convergence_at_some_time}
\end{figure}

For some specific time steps and some specific initial perturbation amplitudes, we observed the non-convergence of all numerical schemes investigated here. For example, \cref{fig:non_convergence_at_some_time} shows a non-convergence of the FSM scheme for the final time $T = 10000$ and initial condition \eqref{eqn:sin02}, with an initial amplitude $\hat{A} = 15.2$, a space step $\Delta x = 2^{-7}$, and a time step $\Delta t = 2^{-7}$. In this case, for $t > 4000$, the error $E(t)$ undergoes oscillations within a narrow range $[\num{2e-4}, \num{8e-4}]$, as shown in sub-panel (a). The inset in sub-panel (a) presents the zoom-in of $E(t)$ in the time interval $[7000, 7300]$, to help us better visualize the variation in this error function. This non-convergence persists if we run simulations for longer time, such as $T = 50000$ or $T = 100000$.  \cref{fig:non_convergence_at_some_time}(b), (c), and (d) shows total densities at times $t = 4300$, $t = 7000$, and $t = 8500$, respectively. These total densities have different values of $E(t)$ in range $[\num{2e-4}, \num{8e-4}]$, but all patterns are non-symmetric.

\begin{table}[!ht]
        \centering
	\caption{Examples of non-convergence of numerical schemes. These simulations were performed using the initial condition \eqref{eqn:sin02}, with an initial amplitude of $\hat{A} = 25.0$ and a space step $\Delta x = 2^{-7}$. For the cases where the solution converges, we also specify the symmetry of this solution, i.e., ``EVEN" symmetry. }
	\label{table:non_convergence}
	\begin{tabular}{|l|c|c|c|}
		\hline
		Scheme used & $\Delta t = 2^{-7}$ & $\Delta t = 2^{-6}$ & $\Delta t = 2^{-5}$ \\ \hline
		Upwind & Non-convergence  & Non-convergence & Non-convergence	\\ \hline
		MacCormack & EVEN  & EVEN  & Non-convergence	\\ \hline
		FSM & Non-convergence  & Non-convergence & Non-convergence	\\ \hline
		QSA & Non-convergence  & Non-convergence & Non-convergence	\\ \hline
		QSA\_Center & EVEN  & EVEN  & Non-convergence	\\ \hline
		QSA\_BW & Non-convergence & Non-convergence & EVEN 	\\ \hline
		QSA\_LW & EVEN  & Non-convergence & Non-convergence	\\ \hline
		QSA\_Minmod & EVEN  & EVEN & Non-convergence	\\ \hline
		QSA\_Superbee & Non-convergence & Non-convergence & EVEN	\\ \hline
		QSA\_MC & EVEN  & Non-convergence & EVEN 	\\ \hline
	\end{tabular}
\end{table}

The non-convergence occurs for all numerical schemes when the initial perturbation amplitude is large enough, typically $\hat{A}>10$. However, this threshold depends on both the time step and the numerical scheme employed; see \cref{fig:all_schemes_7_5_6} discussed in more detail in the next subsection. In \cref{table:non_convergence}, we summarise the results of several simulations performed using an initial amplitude $\hat{A}=25.0$, a spatial step $\Delta x=2^{-7}$, and time steps $\Delta t={2^{-7}; \, 2^{-6}; \, 2^{-5}}$. For the upwind, FSM, and QSA schemes, non-convergence occurs at all three time steps. For the QSA\_BW, and QSA\_Superbee schemes, non-convergence occurs at two smaller time steps $\Delta t=2^{-7}$ and $\Delta t=2^{-6}$, while for the time step $\Delta t=2^{-5}$ the solution converges to an even-symmetric localised pattern. In contrast, the QSA\_LW scheme converges to an even symmetric solution for the time step $\Delta t=2^{-7}$, while the non-convergence occurs at two larger time steps $\Delta t=2^{-6}$ and $\Delta t=2^{-5}$.
For the MacCormack, the QSA\_Center, the QSA\_Minmod schemes, convergence occurs only at time step $\Delta t=2^{-5}$, while the QSA\_MC scheme is non-convergent at time step $\Delta t=2^{-6}$.
Therefore, the non-convergence of the numerical schemes is somewhat random (i.e., we could not identify the exact rules). Increasing or decreasing the time steps or the initial amplitude cannot control the non-convergence, as shown in \cref{table:non_convergence} and \cref{fig:all_schemes_7_5_6}.


\subsection{Change in the solution symmetry as we vary the amplitude of initial conditions, for different numerical schemes. }\label{subsec:overview}

\begin{figure}[!ht]
        \centering
        \begin{adjustbox}{clip,trim=1.3cm 1.0cm 1.2cm 1.5cm}
	\resizebox{1.0\textwidth}{!}{\input{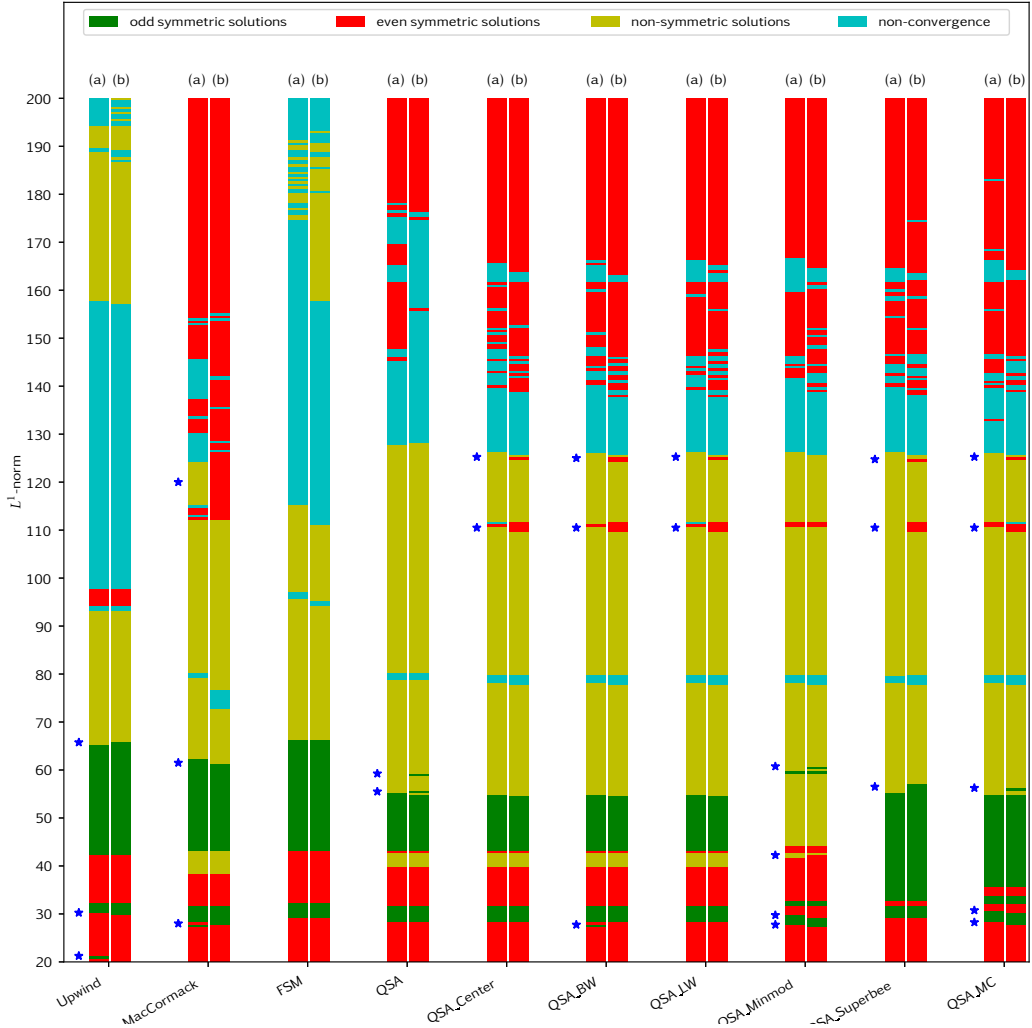}}
        \end{adjustbox}
        \caption{The change of types of solutions for each scheme as moving along the $\mathit{L}^1$-norm of total density. A total of 7220 simulations were executed by using the initial condition \eqref{eqn:sin02} with a space step of $\Delta x = 2^{-7}$ and with two different time steps of $\Delta t = 2^{-5}$, and  $\Delta t = 2^{-6}$, see columns labeled as \textbf{(a)} and \textbf{(b)}, respectively. For each scheme and each time step, 361 simulations were executed with 360 different values of the initial amplitude $\hat{A} \in K = \{0.001; \, 0.1 \times n\}$ with $n = \overline{1,36}$, corresponding to $\normLp{u}{1} \in (20, 200]$, see \eqref{eqn:L1-A}. To show the continuous $\mathit{L}^1$-norm in this figure, we assume that the solution obtained with the initial amplitude of $\hat{A} = k \in K$ represents all solutions with initial amplitudes inside the interval $ \hat{A} \in (k - 0.05, k+0.05]$. The blue stars on the left side of each scheme indicate positions where different solution types (odd symmetric, even symmetric, or non-symmetric) occur between two time step values.} \label{fig:all_schemes_7_5_6}
\end{figure}

To conclude this section, we now investigate the changes in solution symmetries when the $\mathit{L}^1$-norm of the solutions is increased from 20 to 200. To do this, for each scheme with a fixed space-time step, we perform 361 simulations with different initial amplitudes which belong to the set $K = \{0.001; \, 0.1 \times n\}$ where $n = \overline{1,36}$. From \eqref{eqn:L1-A}, we obtain the $\mathit{L}^1$-norm of solutions  belongs to $\{20.005; \, 20 + 5 \times n\} \subset (20, 200]$.  Further, to show the continuous $\mathit{L}^1$-norm, we assume that the solution obtained with the initial amplitude of $\hat{A} = k \in K$ is actually an average of solutions with initial amplitudes in the interval $ \hat{A} \in (k - 0.05, k+0.05]$. This assumption usually holds when the simulation converges to a solution. However, it does not really hold in the non-convergence case: i.e., we could have non-convergence for $\hat{A} = k$, convergence for $\hat{A} = k + 0.01$ and non-convergence again for $\hat{A} = k + 0.02$. But since we can only execute a finite and small number of test cases, we need this assumption to have a continuous $\mathit{L}^1$-norm.

\cref{fig:all_schemes_7_5_6}  presents the change in the solution symmetries for space step $\Delta x = 2^{-7}$ and two different time steps $\Delta t = 2^{-5}$ and $\Delta t = 2^{-6}$; see columns labeled as \textbf{(a)} and \textbf{(b)}, respectively. The blue stars on the left side of each scheme in \cref{fig:all_schemes_7_5_6} indicate positions where different solution types (odd symmetric, even symmetric, or non-symmetric) occur between two different time step values. In view of creating (as future work) a bifurcation diagram for the transitions between solutions with different types of symmetries, we note that the transitions between the odd-symmetric (green colour) and even-symmetric (red colour) or non-symmetric (yellow colour) solutions occur for different amplitudes of initial perturbations, and are different for different numerical schemes. This suggests that the corresponding bifurcation diagrams would be slightly different for different schemes.

Regarding the non-convergence of numerical schemes (light blue colour in \cref{fig:all_schemes_7_5_6}), the MacCormack scheme exhibits the smallest $\mathit{L}^1$-norm of $73$, and interestingly, it has the fewest instances of non-convergence. On the other hand, the upwind and FSM schemes experience non-convergence at a later stage (around an $\mathit{L}^1$-norm of 94), and they have the highest number of non-convergence cases. In terms of solution types for large $\mathit{L}^1$-norm (ranging from $140$ to $200$), the upwind and FSM schemes yield non-symmetric solutions, while the other schemes produce even symmetric solutions. Additionally, when the $\mathit{L}^1$-norm ranges from 20 to 140, noticeable differences are observed between the various numerical schemes and their corresponding time steps. 
Despite these differences, the group of QSA, QSA\_Center, QSA\_BW, and QSA\_LW schemes exhibit the smallest discrepancy. Particularly, QSA\_Center and QSA\_LW demonstrate a very slight distinction between the two time steps of $\Delta t = 2^{-5}$, and  $\Delta t = 2^{-6}$. However, in our simulations, the QSA\_Center scheme has the advantage of faster convergence compared to the QSA\_LW scheme. For more detail, please refer to \cref{fig:all_schemes_7_5_6}.

%% file: Conclusion.tex
\section{Summary and Conclusions} \label{sec:conclusion}

In this paper, we studied the effect of different numerical schemes (and their time and space steps), as well as the effect of different initial conditions, on the transitions between different localised solutions characterized by different symmetries, which were previously observed in non-local hyperbolic systems for ecological aggregations~\cite{Eftimie-Vries-Lewis-2009}. We have emphasised here several numerical issues. First, we observed the presence of two distinct types of numerical solutions (transient and steady-state solutions) that exhibited very small errors and could be misleading in terms of stopping times for continuation algorithms (which should stop at steady-state solutions and not transient solutions). This also implies an extremely slow convergence.  Second, in some cases, none of the investigated numerical schemes converged, posing a challenge to the numerical analysis. Lastly, we have discovered that the choice of numerical schemes (and the time and space steps used for simulations) as well as the choice of the initial conditions, exert a significant influence on the type and symmetry of bifurcating solutions. In consequence, the resulting bifurcation diagrams (which will be developed in the future) may vary when different numerical schemes and/or corresponding parameters are employed (see \cref{fig:all_schemes_7_5_6}). 

\begin{figure}[!ht]
	\centering
	\resizebox{0.9\textwidth}{!}{\input{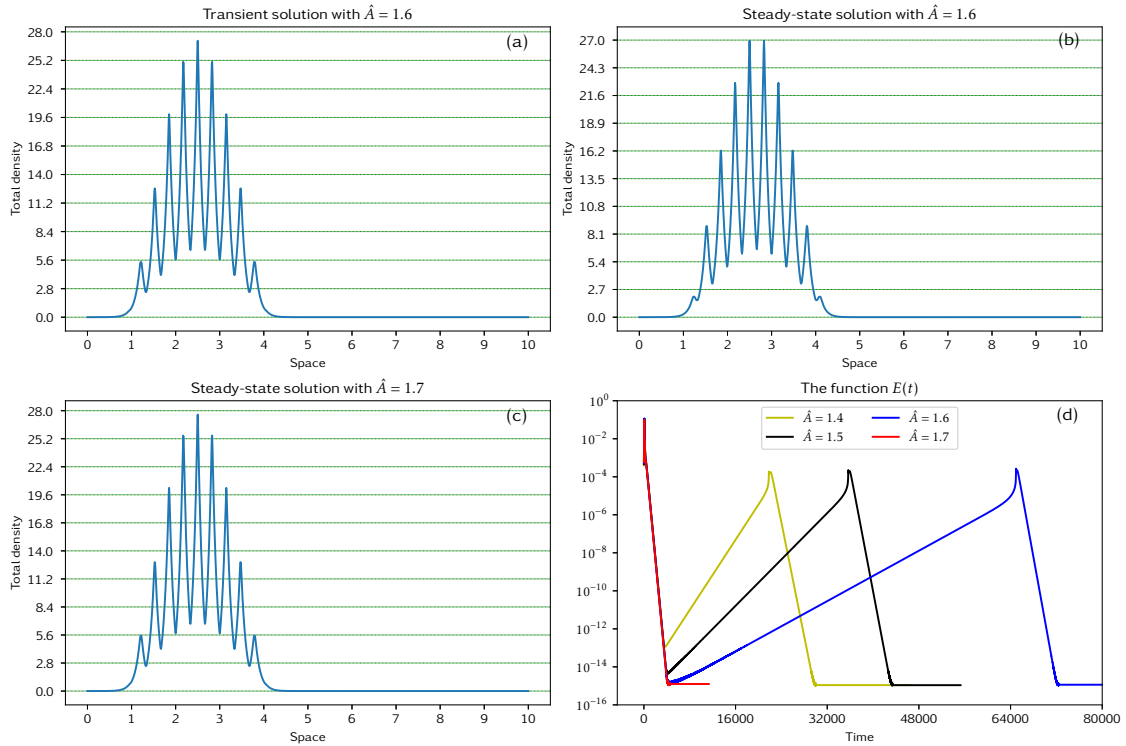}}
        \caption{Convergence history for several initial amplitudes $\hat{A} = 1.4; 1.5; 1.6; 1.7$. These simulations were performed using the QSA\_Center scheme, a space step of $\Delta x = 2^{-7}$ and a time step of $\Delta t = 2^{-5}$.}
	\label{fig:open_question}
\end{figure}

{\bf Open questions.} In the following, we briefly mention some of the open questions identified through this study, which will have to be addressed in the future:
\begin{itemize}
\item We need to investigate in more detail (numerically and theoretically) the slow convergence to steady state solutions, and the differences in the types of solutions when we change the time/space steps.
\item Theoretical investigation of the convergence of all these numerical schemes for nonlocal hyperbolic systems is necessary,  to clarify whether: (i) is it the scheme that is not converging for those specific time and/or space steps, or (ii) is it the dynamics of the system that could be exhibiting another type of solution -- heteroclinic orbits connecting different localised solutions -- in that parameter region, and the numerical approach used does not capture these extra dynamics.
\item The impact of the very small transient error could lead to the stopping of some numerical continuation algorithms, even when the solution did not reach a steady state. Usually, many numerical algorithms stop if the error tolerance (equivalent to $E(t)$) is less than \num{e-8} or \num{e-10}), but as we have seen in this study, even if we chose the tolerance small enough (here \num{e-14}), in some cases we cannot ensure that we really reach a steady-state solution (see also the discussion in Appendix A.1). For instance, for initial amplitudes $\hat{A} = 1.4, \; 1.5, \; 1.6$, we obtain even symmetric steady-state solutions (see \cref{fig:open_question}(b)), and in these cases, we also obtain odd symmetric transient solutions (see \cref{fig:open_question}(a)). We observe that the error of transient solutions is decreasing as we increase the initial amplitudes, as shown by yellow, black, and blue curves in \cref{fig:open_question}(d). However, if the initial amplitude is $\hat{A}=1.7$, we obtain an odd symmetric steady-state solution (see \cref{fig:open_question}(c)) and no transient solution (see red curve in  \cref{fig:open_question}(d)). This raises the question of whether this is a stable steady-state solution or a transient solution; but because of the limitations of the numerical scheme, as well as the \textbf{double} data type in \textbf{C++}, we cannot answer this question.
\end{itemize}

%% file: Appendix.tex
\section*{Appendix A.1. Numerical transient solutions and numerical steady-state solutions} \label{sec:appendix_A1}

In this Appendix, we show a few more numerical results emphasizing the issue of having transient vs. steady-state solutions with different symmetries presented in \cref{subsec:transient_solution}.

\begin{figure}[!ht] 
	\centering
        \resizebox{\textwidth}{!}{\input{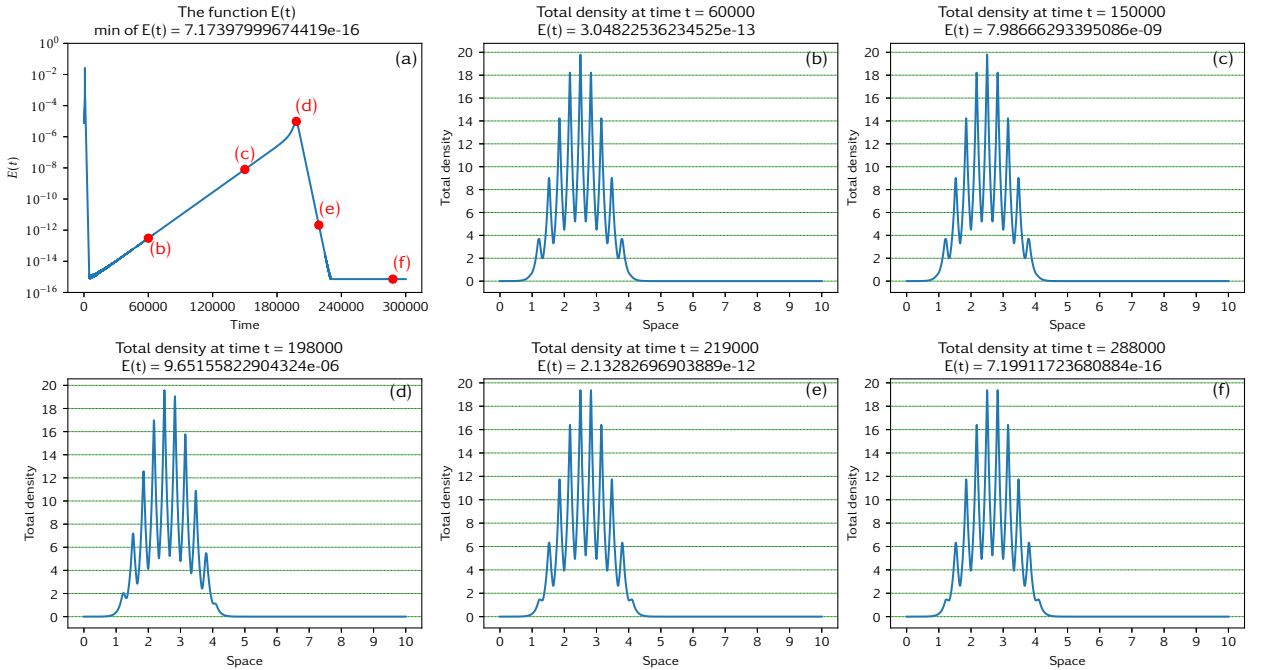}}
	\caption{Total densities at the varying time along the error $E(t)$ in the case where the transient solution exists. The simulation used the upwind scheme and the initial condition \eqref{eqn:sin02} with an initial amplitude of $\hat{A} = 0.2$, a space step of $\Delta x = 2^{-7}$ and a time step of $\Delta t = 2^{-6}$. }
	\label{fig:transient_at_some_time}
\end{figure}

As shown in \cref{fig:transient_steady_state_solution}, there is an odd symmetric transient solution at time $t = 6000$ with $E(t) = \num{1.11 e-15}$ and an even symmetric steady-state solution at time $t = 240 000$ with $E(t) = \num{7.2 e-16}$. To have a better understanding of what happens in this case, in \cref{fig:transient_at_some_time} we show the types of solutions observed at different times; see sub-panel (a). Sub-panel (b) shows an odd symmetric solution at time $t = 60 000$ with $E(t) = \num{3.04 e-13}$. At time $t = 150 000$ the total density is still odd symmetric with $E(t) = \num{7.99 e-9}$; see sub-panel (c). At the local maximum of $E(t)$, we obtain a non-symmetric total density solution (at time $t = 198 000$) with $E(t) = \num{9.65 e-6}$; see sub-panel (d). Note that for all simulations performed in this paper, the range of the local maximum of $E(t)$ is between $10^{-6}$ and $10^{-2}$. An even symmetric total density is observed in sub-panel (e) at time $t = 219 000$ with $E(t) = \num{2.13 e-12}$. From time $t > 220 000$, we obtain only even-symmetric steady-state solution with $E(t) = \num{7.2 e-16}$, as shown in sub-panel (f).

\begin{figure}[!ht] 
	\centering
        \resizebox{0.9\textwidth}{!}{\input{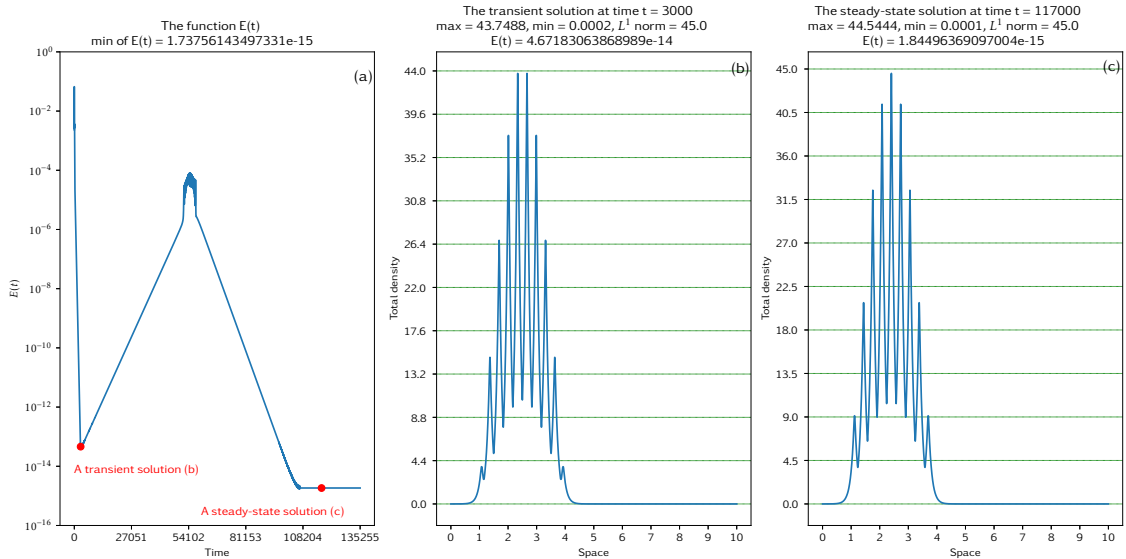}}
	\caption{An even symmetric transient solution and an odd symmetric steady-state solution. The simulation used the QSA\_Minmod scheme and the initial condition \eqref{eqn:sin02} with an initial amplitude of $\hat{A} = 5.0$, a space step of $\Delta x = 2^{-9}$ and a time step of $\Delta t = 2^{-7}$. }
	\label{fig:transient_steady_state_solution_QSA_Minmod_5}
\end{figure}

\begin{figure}[!ht] 
	\centering
        \resizebox{0.9\textwidth}{!}{\input{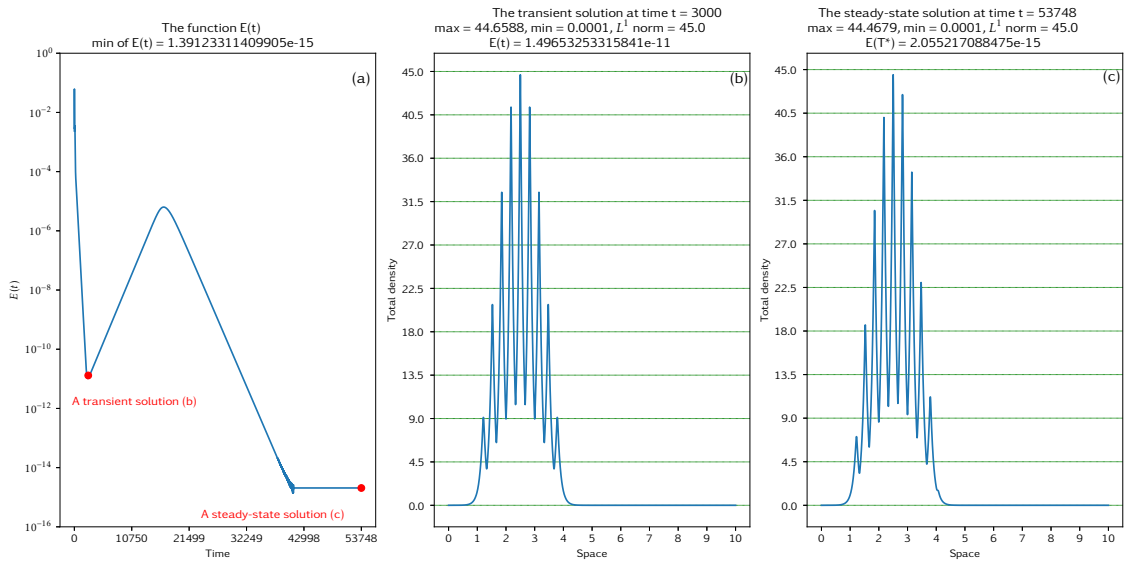}}
	\caption{An odd symmetric transient solution and a non-symmetric steady-state solution. The simulation used the QSA scheme and the initial condition \eqref{eqn:sin02} with an initial amplitude of $\hat{A} = 5.0$, a space step of $\Delta x = 2^{-8}$ and a time step of $\Delta t = 2^{-7}$. }
	\label{fig:transient_steady_state_solution_QSA_4}
\end{figure}

Moreover, in \cref{fig:transient_steady_state_solution} we observed a transition from odd-symmetric to even-symmetric solutions; in contrast, in \cref{fig:transient_steady_state_solution_QSA_Minmod_5} (a simulation of the QSA\_Minmod scheme with an initial amplitude of $\hat{A} = 5.0$, a space step of $\Delta x = 2^{-9}$ and a time step $\Delta t = 2^{-7}$) we can show the opposite transition: from even-symmetric to odd-symmetric solutions. In some cases we could also observe a transition from either odd or even symmetric to non-symmetric solutions, or even the opposite transition; see, for example, \cref{fig:transient_steady_state_solution_QSA_4} (a simulation of the QSA scheme with an initial amplitude of $\hat{A} = 5.0$, a space step of $\Delta x = 2^{-8}$ and a time step $\Delta t = 2^{-7}$). \\
Overall, transient solutions were observed to exist for all schemes analyzed in this paper, with varying initial amplitudes and different time and space steps; for example, see \cref{fig:transient_steady_state_solution},  \cref{fig:test_dx_error}, \cref{fig:0300_some_schemes}, \cref{fig:open_question}, \cref{fig:transient_at_some_time}, \cref{fig:transient_steady_state_solution_QSA_Minmod_5}, and \cref{fig:transient_steady_state_solution_QSA_4}.

\section*{Appendix A.2. The influence of the time step} \label{sec:appendix_A2}

In this Appendix, we show a few more numerical results emphasizing the impact of the time step on different symmetries solutions presented in \cref{subsec:time_step}.

\begin{figure}[!ht] 
	\centering
        \resizebox{0.9\textwidth}{!}{\input{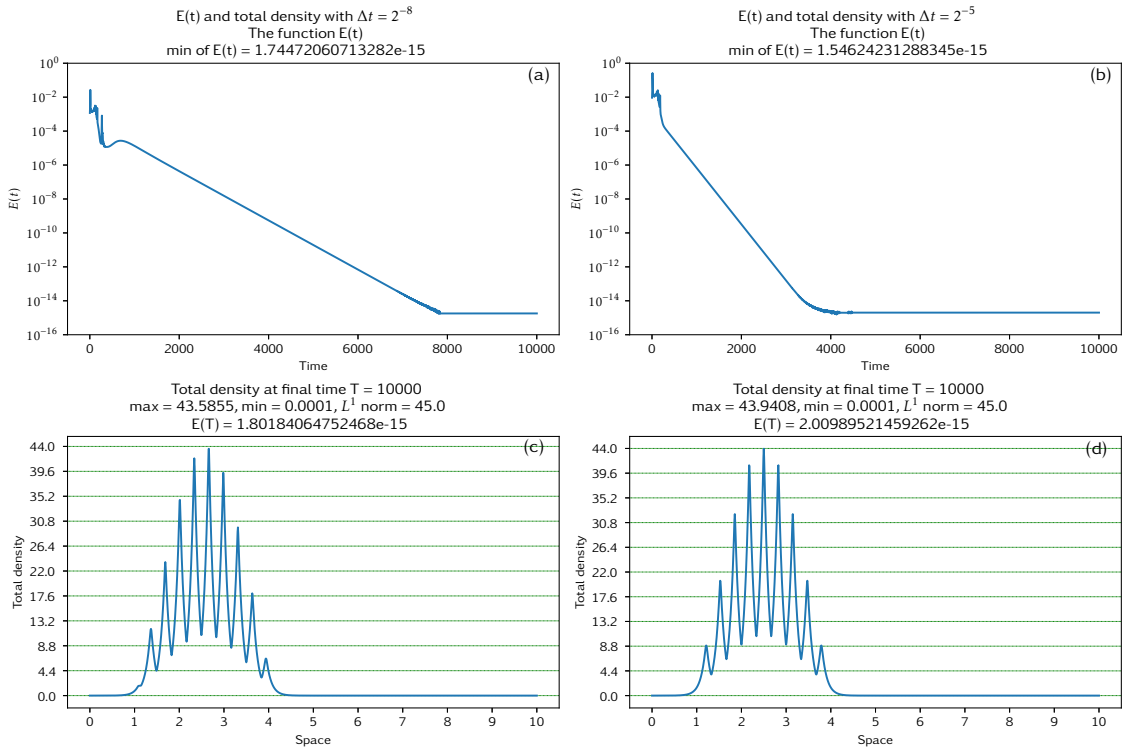}}
	\caption{In the first column, a non-symmetric solution is obtained with a time step of $\Delta t = 2^{-8}$, however, in the second column, an odd symmetric solution is obtained by using a time step of $\Delta t = 2^{-5}$. These simulations were executed by using the QSA\_MC scheme, and the initial condition \eqref{eqn:sin02} with an initial amplitude of $\hat{A} = 5.0$, and a space step of $\Delta x = 2^{-7}$.}
	\label{fig:test_dt_5000_QSA_MC}
\end{figure}

In \cref{fig:test_dt_0200_upwind}, we have observed differences in the symmetry of localised solution when varying the time steps: an odd-symmetric steady-state solution for $\Delta t = 2^{-5}$, and an even-symmetric steady-state solution for $\Delta t = 2^{-6}$. Note that, for $\Delta t = 2^{-6}$, there exists also an odd symmetric transient solution. However, it is important to note that differences in the type of solutions generated by varying time steps exist also in cases where there is no transient solution.
For example, using the QSA\_MC scheme with an initial amplitude of $\hat{A} = 5.0$, time steps of $\Delta t = 2^{-8}$ and $\Delta t = 2^{-5}$ produce a non-symmetric solution as shown in \cref{fig:test_dt_5000_QSA_MC} (c), and an odd-symmetric solution as shown in \cref{fig:test_dt_5000_QSA_MC}(d). In these two simulations, there is no transient solution, see \cref{fig:test_dt_5000_QSA_MC}(a) and (b).

\begin{figure}[!ht] 
	\centering
        \resizebox{0.9\textwidth}{!}{\input{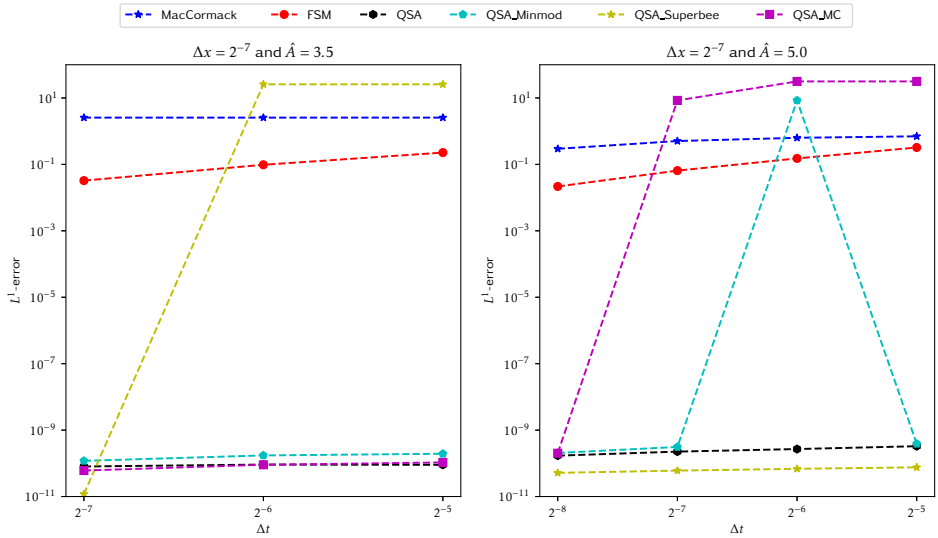}}
	\caption{Comparison between solutions at several time steps and the numerical reference solution for some schemes. These simulations were executed by using the initial condition \eqref{eqn:sin02} with initial amplitudes of $\hat{A} = 3.5$ and $\hat{A} = 5.0$, a space step of $\Delta x = 2^{-7}$.}
	\label{fig:test_dt_error}
\end{figure}
\begin{figure}[!ht] 
	\centering
        \resizebox{0.85\textwidth}{!}{\input{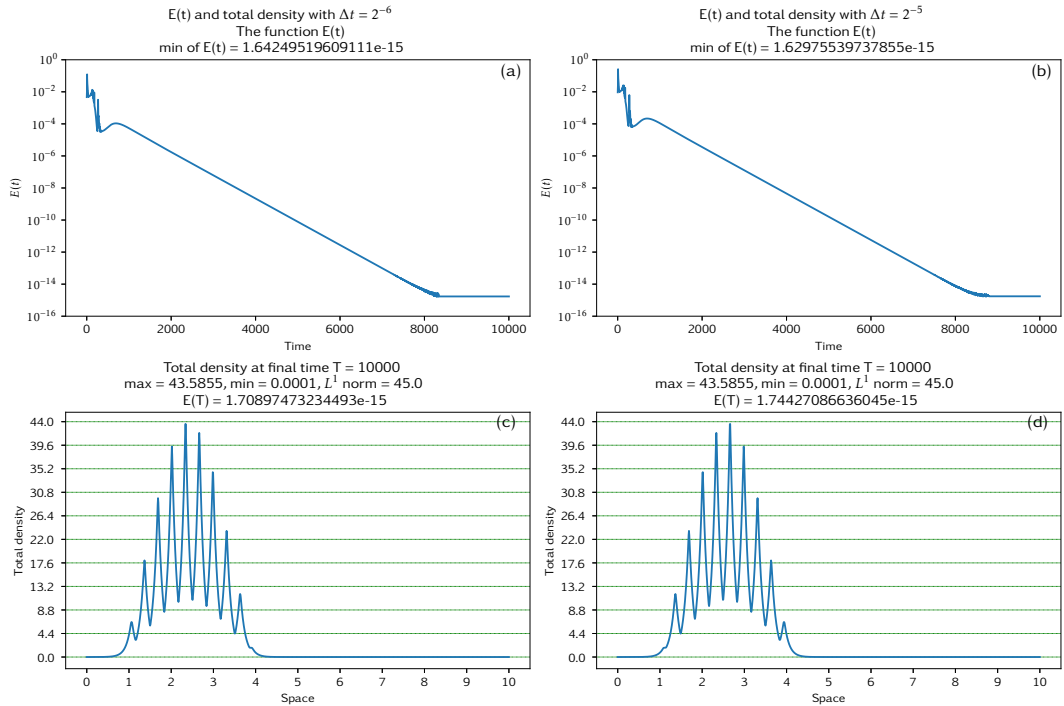}}
	\caption{In both cases $\Delta t = 2^{-6}$ and $\Delta t = 2^{-5}$, non-symmetric solutions are obtained, however, there is a slight difference in the geometric of the two solutions. These simulations were executed by using the QSA\_Minmod scheme, and the initial condition \eqref{eqn:sin02} with an initial amplitude of $\hat{A} = 5.0$, and a space step of $\Delta x = 2^{-7}$.}
	\label{fig:test_dt_5000_QSA_Minmod}
\end{figure}

To compare solutions at different time steps, we have chosen two cases as our reference numerical solutions: (i) $\Delta t = 2^{-8}$ and initial amplitude $\hat{A} = 3.5$; and (ii) $\Delta t = 2^{-9}$ and initial amplitude of $\hat{A} = 5.0$. \cref{fig:test_dt_error} displays the $\mathit{L}^1$-norm of the difference between these reference solutions and the solutions obtained with other time steps.
Although the MacCormack and FSM schemes produce similar solution types, their $\mathit{L}^1$-errors are relatively large, typically greater than $\num{e-2}$ (this is probably caused by the fact that they are two-stage schemes). For the QSA\_Superbee scheme with the initial amplitude  $\hat{A} = 3.5$ and the QSA\_MC scheme with the initial amplitude $\hat{A} = 5.0$, the $\mathit{L}^1$-errors are also large due to changes in the solution types. On the other hand, schemes such as upwind, QSA, QSA\_Center, QSA\_BW, and QSA\_LW produce similar solutions types with very small $\mathit{L}^1$-errors, around $\num{e-10}$. For brevity, in \cref{fig:test_dt_error} we only display the error of the QSA scheme.
Finally, in \cref{fig:test_dt_5000_QSA_Minmod} we show that for the QSA\_Minmod scheme with initial amplitude $\hat{A} = 5.0$, all solutions are non-symmetric, but with geometric differences between the peaks of the solutions, resulting in a large error.

%% file: main.bbl
\begin{thebibliography}{10}

\bibitem{Aldana-Dossetti-Huepe-Kenkre-Larralde-2007}
M.~Aldana, V.~Dossetti, C.~Huepe, V.~M. Kenkre, and H.~Larralde.
\newblock Phase transitions in systems of self-propelled agents and related
  network models.
\newblock {\em Phys. Rev. Lett.}, 98:095702, Mar 2007.

\bibitem{Avitabile-Lloyd-Burke-Knobloch-Sandstede-2010}
Daniele Avitabile, David J.~B. Lloyd, John Burke, Edgar Knobloch, and Bj\"{o}rn
  Sandstede.
\newblock To snake or not to snake in the planar {S}wift-{H}ohenberg equation.
\newblock {\em SIAM J. Appl. Dyn. Syst.}, 9(3):704--733, 2010.

\bibitem{Bale-Leveque-Mitran-Rossmanith-2002}
Derek~S. Bale, Randall~J. Leveque, Sorin Mitran, and James~A. Rossmanith.
\newblock A wave propagation method for conservation laws and balance laws with
  spatially varying flux functions.
\newblock {\em SIAM J. Sci. Comput.}, 24(3):955--978, 2002.

\bibitem{Beck-Knobloch-Lloyd-Sandstede-Wagenknecht-2009}
Margaret Beck, J\"{u}rgen Knobloch, David J.~B. Lloyd, Bj\"{o}rn Sandstede, and
  Thomas Wagenknecht.
\newblock Snakes, ladders, and isolas of localized patterns.
\newblock {\em SIAM J. Math. Anal.}, 41(3):936--972, 2009.

\bibitem{Buono-Eftimie-2014-MMMA}
P.-L. Buono and R.~Eftimie.
\newblock Analysis of {H}opf/{H}opf bifurcations in nonlocal hyperbolic models
  for self-organised aggregations.
\newblock {\em Math. Models Methods Appl. Sci.}, 24(2):327--357, 2014.

\bibitem{Buono-Eftimie-2014-SIAM}
Pietro-Luciano Buono and Raluca Eftimie.
\newblock Codimension-two bifurcations in animal aggregation models with
  symmetry.
\newblock {\em SIAM J. Appl. Dyn. Syst.}, 13(4):1542--1582, 2014.

\bibitem{Burke-Knobloch-2006}
John Burke and Edgar Knobloch.
\newblock Localized states in the generalized {S}wift-{H}ohenberg equation.
\newblock {\em Phys. Rev. E (3)}, 73(5):056211, 15, 2006.

\bibitem{Burke-Knobloch-2007-Chaos}
John Burke and Edgar Knobloch.
\newblock Homoclinic snaking: structure and stability.
\newblock {\em Chaos}, 17(3):037102, 15, 2007.

\bibitem{Burke-Knobloch-2007-PLA}
John Burke and Edgar Knobloch.
\newblock Snakes and ladders: localized states in the {S}wift-{H}ohenberg
  equation.
\newblock {\em Phys. Lett. A}, 360(6):681--688, 2007.

\bibitem{Chuang-Orsogna-Marthaler-Bertozzi-Chayes-2007}
Yao-li Chuang, Maria~R. D'Orsogna, Daniel Marthaler, Andrea~L. Bertozzi, and
  Lincoln~S. Chayes.
\newblock State transitions and the continuum limit for a 2{D} interacting,
  self-propelled particle system.
\newblock {\em Phys. D}, 232(1):33--47, 2007.

\bibitem{PyDSTool}
R.H. Clewley, W.E. Sherwood, M.D. LaMar, and J.~Guckenheimer.
\newblock Pydstool: a software environment for dynamical systems modeling.
\newblock http://pydstool.sourceforge.net, 2010.

\bibitem{Couzin-Krause-James-Ruxtion-2002}
Iain~D. Couzin, Jens Krause, Richard James, Graeme~D. Ruxton, and Nigel~R.
  Franks.
\newblock Collective memory and spatial sorting in animal groups.
\newblock {\em J. Theoret. Biol.}, 218(1):1--11, 2002.

\bibitem{Czirok-Barabasi-Vicsek-1999}
Andr\'as Czir\'ok, Albert-L\'aszl\'o Barab\'asi, and Tam\'as Vicsek.
\newblock Collective motion of self-propelled particles: Kinetic phase
  transition in one dimension.
\newblock {\em Phys. Rev. Lett.}, 82:209--212, Jan 1999.

\bibitem{Czirok-Stanley-Vicsek-1997}
András Czirók, H~Eugene Stanley, and Tamás Vicsek.
\newblock Spontaneously ordered motion of self-propelled particles.
\newblock {\em Journal of Physics A: Mathematical and General}, 30(5):1375, mar
  1997.

\bibitem{MATCONT}
A.~Dhooge, W.~Govaerts, and Yu.~A. Kuznetsov.
\newblock M{ATCONT}: a {MATLAB} package for numerical bifurcation analysis of
  {ODE}s.
\newblock {\em ACM Trans. Math. Software}, 29(2):141--164, 2003.

\bibitem{AUTO-07p}
E.J. Doedel, A.~Champneys, F.~Dercole, T.~Fairgrieve, Y.~Kuznetsov, B.~Oldeman,
  R.~Paffenroth, B.~Sandstede, X.~Wang, and C.~Zhang.
\newblock Auto 2007p: Continuation and bifurcation software for ordinary
  differential equations (with homcont).
\newblock http://cmvl.cs.concordia.ca/auto, 2007.

\bibitem{Eftimie-Vries-Lewis-2007}
R.~Eftimie, G.~de~Vries, and M.~A. Lewis.
\newblock Complex spatial group patterns result from different animal
  communication mechanisms.
\newblock {\em Proc. Natl. Acad. Sci. USA}, 104(17):6974--6979, 2007.

\bibitem{Eftimie-Vries-Lewis-2009}
R.~Eftimie, G.~de~Vries, and M.~A. Lewis.
\newblock Weakly nonlinear analysis of a hyperbolic model for animal group
  formation.
\newblock {\em J. Math. Biol.}, 59(1):37--74, 2009.

\bibitem{Eftimie-Vries-Lewis-Lutscher-2007}
R.~Eftimie, G.~de~Vries, M.~A. Lewis, and F.~Lutscher.
\newblock Modeling group formation and activity patterns in self-organizing
  collectives of individuals.
\newblock {\em Bull. Math. Biol.}, 69(5):1537--1565, 2007.

\bibitem{Eftimie-2012}
Raluca Eftimie.
\newblock Hyperbolic and kinetic models for self-organized biological
  aggregations and movement: a brief review.
\newblock {\em J. Math. Biol.}, 65(1):35--75, 2012.

\bibitem{XXP}
B.~Ermentrout.
\newblock Xppaut.
\newblock http://www.math.pitt.edu/~bard/xpp/xpp.html, 2008.

\bibitem{Feder-2007}
Toni Feder.
\newblock {Statistical physics is for the birds}.
\newblock {\em Physics Today}, 60(10):28--30, 10 2007.

\bibitem{Fetecau-2011}
Razvan~C. Fetecau.
\newblock Collective behavior of biological aggregations in two dimensions: a
  nonlocal kinetic model.
\newblock {\em Math. Models Methods Appl. Sci.}, 21(7):1539--1569, 2011.

\bibitem{Fetecau-Eftimie-2010}
Razvan~C. Fetecau and Raluca Eftimie.
\newblock An investigation of a nonlocal hyperbolic model for self-organization
  of biological groups.
\newblock {\em J. Math. Biol.}, 61(4):545--579, 2010.

\bibitem{Flierl-Grunbaum-Levin-Olson-1999}
G.~Flierl, D.~Grünbaum, S.~Levin, and D.~Olson.
\newblock From individuals to aggregations: the interplay between behavior and
  physics.
\newblock {\em Journal of Theoretical Biology}, 196(4):397--454, 1999.

\bibitem{Gueron-Levin-Rubenstein-1996}
Shay Gueron, Simon~A. Levin, and Daniel~I. Rubenstein.
\newblock The dynamics of herds: From individuals to aggregations.
\newblock {\em Journal of Theoretical Biology}, 182(1):85--98, 1996.

\bibitem{Helbing-Treiber-1999}
D.~Helbing and M.~Treiber.
\newblock Numerical simulation of macroscopic traffic equations.
\newblock {\em {Computing in Science \& Engineering}}, 1(5):89--98, 1999.

\bibitem{Kuehn-2015}
Christian Kuehn.
\newblock Efficient gluing of numerical continuation and a multiple solution
  method for elliptic {PDE}s.
\newblock {\em Appl. Math. Comput.}, 266:656--674, 2015.

\bibitem{LeVeque-1998}
Randall~J. LeVeque.
\newblock Balancing source terms and flux gradients in high-resolution
  {G}odunov methods: the quasi-steady wave-propagation algorithm.
\newblock {\em J. Comput. Phys.}, 146(1):346--365, 1998.

\bibitem{LeVeque-2002-Finite-volume}
Randall~J. LeVeque.
\newblock {\em Finite volume methods for hyperbolic problems}.
\newblock Cambridge Texts in Applied Mathematics. Cambridge University Press,
  Cambridge, 2002.

\bibitem{LeVeque-2007-Finite-difference}
Randall~J. LeVeque.
\newblock {\em Finite difference methods for ordinary and partial differential
  equations}.
\newblock Society for Industrial and Applied Mathematics (SIAM), Philadelphia,
  PA, 2007.
\newblock Steady-state and time-dependent problems.

\bibitem{Liu-Xu-2017}
Yongli Liu and Yancong Xu.
\newblock Localized patterns of the swift-{H}ohenberg equation with a
  dissipative term.
\newblock {\em Ann. Appl. Math.}, 33(1):6--17, 2017.

\bibitem{Lukeman-Li-Keshet-2009}
Ryan Lukeman, Yue-Xian Li, and Leah Edelstein-Keshet.
\newblock A conceptual model for milling formations in biological aggregates.
\newblock {\em Bull. Math. Biol.}, 71(2):352--382, 2009.

\bibitem{Lutscher-2002}
Frithjof Lutscher.
\newblock Modeling alignment and movement of animals and cells.
\newblock {\em J. Math. Biol.}, 45(3):234--260, 2002.

\bibitem{MacCormack-2003}
Robert~W. MacCormack.
\newblock The effect of viscosity in hypervelocity impact cratering.
\newblock {\em Journal of Spacecraft and Rockets}, 40(5):757--763, 2003.

\bibitem{Mogilner-Keshet-1999}
Alexander Mogilner and Leah Edelstein-Keshet.
\newblock A non-local model for a swarm.
\newblock {\em J. Math. Biol.}, 38(6):534--570, 1999.

\bibitem{Okubo-Grunbaum-Keshet-2001}
Akira Okubo, Daniel Gr{\"u}nbaum, and Leah Edelstein-Keshet.
\newblock {\em The Dynamics of Animal Grouping}, pages 197--237.
\newblock Springer New York, New York, NY, 2001.

\bibitem{Othmer-Dunbar-Alt-1988}
H.~G. Othmer, S.~R. Dunbar, and W.~Alt.
\newblock Models of dispersal in biological systems.
\newblock {\em J. Math. Biol.}, 26(3):263--298, 1988.

\bibitem{Parrish-Keshet-1999}
Julia~K. Parrish and Leah Edelstein-Keshet.
\newblock Complexity, pattern, and evolutionary trade-offs in animal
  aggregation.
\newblock {\em Science}, 284(5411):99--101, 1999.

\bibitem{Pfistner-1990}
Beate Pfistner.
\newblock {\em A One Dimensional Model for the Swarming Behavior of
  Myxobacteria}, pages 556--565.
\newblock Springer Berlin Heidelberg, Berlin, Heidelberg, 1990.

\bibitem{Reynolds-1987}
Craig~W. Reynolds.
\newblock Flocks, herds and schools: A distributed behavioral model.
\newblock {\em SIGGRAPH Comput. Graph.}, 21(4):25–34, aug 1987.

\bibitem{Schmidt-Avitabile-2020}
Helmut Schmidt and Daniele Avitabile.
\newblock Bumps and oscillons in networks of spiking neurons.
\newblock {\em Chaos}, 30(3):033133, 13, 2020.

\bibitem{Taylor-Dawes-2010}
Chris Taylor and Jonathan~H.P. Dawes.
\newblock Snaking and isolas of localised states in bistable discrete lattices.
\newblock {\em Physics Letters A}, 375(1):14--22, 2010.

\bibitem{Topaz-Bertozzi-Lewis-2006}
Chad~M. Topaz, Andrea~L. Bertozzi, and Mark~A. Lewis.
\newblock A nonlocal continuum model for biological aggregation.
\newblock {\em Bull. Math. Biol.}, 68(7):1601--1623, 2006.

\bibitem{Uecker-2022}
Hannes Uecker.
\newblock Continuation and bifurcation in nonlinear {PDE}s---algorithms,
  applications, and experiments.
\newblock {\em Jahresber. Dtsch. Math.-Ver.}, 124(1):43--80, 2022.
\newblock [Vol. 124 (2021) on first page].

\bibitem{Vicsek-Czirok-Jacob-Cohen-Shochet-1995}
Tam\'as Vicsek, Andr\'as Czir\'ok, Eshel Ben-Jacob, Inon Cohen, and Ofer
  Shochet.
\newblock Novel type of phase transition in a system of self-driven particles.
\newblock {\em Phys. Rev. Lett.}, 75:1226--1229, Aug 1995.

\end{thebibliography}
